\documentclass[11pt]{article}

\usepackage{geometry}
\usepackage[utf8]{inputenc}
\usepackage{amsthm, amssymb, natbib, graphicx}
\usepackage{listings}
\usepackage[ruled,vlined]{algorithm2e}
\usepackage{caption}
\usepackage{amsmath,tikz}
\usepackage{array, makecell}
\usepackage{sidecap, caption}
\usepackage{xcolor, color} 
\usepackage{comment} 
\usepackage{algorithm2e}
\usepackage[most]{tcolorbox}
\setcitestyle{numbers, comma,square}
\usepackage{bbm}
\usepackage{caption}
\usepackage{tabularray}
\usepackage{diagbox}
\UseTblrLibrary{booktabs}

\newtheorem{theorem}{Theorem}[section]
\newtheorem{corollary}{Corollary}[section]
\newtheorem{lemma}[theorem]{Lemma}

\newtheorem{remark}{Remark}[section]
\newtheorem{assumption}{Assumption}

\DeclareMathOperator*{\argmin}{argmin}

\DeclareMathOperator{\sign}{sign}

\newcommand{\R}{\mathbb{R}}

\newcommand{\rr}[1]{\textcolor{red}{#1}}

\usepackage[normalem]{ulem}
\newcommand{\bb}[1]{\textcolor{blue}{#1}}

\newcommand{\ignore}[1]{}

\title{Cubic-quartic regularization models for solving polynomial subproblems in third-order tensor methods}

\author{\thanks{wenqi.zhu@maths.ox.ac.uk, University of Oxford, UK} Wenqi Zhu,  \thanks{cartis@maths.ox.ac.uk, University of Oxford, UK} Coralia Cartis}

\author{Coralia Cartis\thanks{The order of the authors is alphabetical; the second author (Wenqi Zhu) is the primary contributor.} \textsuperscript{\normalfont ,}\thanks{Mathematical Institute, Woodstock Road, University of Oxford, Oxford, UK, OX2 6GG.  {\tt cartis@maths.ox.ac.uk} } \quad and \quad
Wenqi Zhu\footnotemark[1] \textsuperscript{\normalfont,}\thanks{Mathematical Institute, Woodstock Road, University of Oxford, Oxford, UK, OX2 6GG.
{\tt wenqi.zhu@maths.ox.ac.uk}}}

\begin{document}
\maketitle
\begin{abstract}
High-order tensor methods for solving both convex and nonconvex optimization problems have recently generated significant research interest, due in part to the natural way in which higher derivatives can be incorporated into adaptive regularization frameworks, leading to algorithms with optimal global rates of convergence and local rates that are faster than Newton's method.
On each iteration, to find the next solution approximation, these methods require the unconstrained local minimization of a (potentially nonconvex) multivariate polynomial of degree higher than two, constructed using third-order (or higher) derivative information, and regularized by an appropriate power of the change in the iterates. Developing efficient techniques for the solution of such subproblems is currently, an ongoing topic of research,  and this paper  addresses this question for the case of the third-order tensor subproblem. In particular, 
we propose the CQR algorithmic framework, for minimizing a nonconvex Cubic multivariate polynomial with  Quartic Regularisation, by sequentially minimizing a sequence of local quadratic models that also incorporate both simple cubic and quartic terms. 
The role of the cubic term is to crudely approximate local tensor information, while the quartic one provides model regularization and controls progress. We provide necessary and sufficient optimality conditions that fully characterise the global minimizers of these cubic-quartic models. We then turn these conditions into secular equations that can be solved using nonlinear eigenvalue techniques. We  show, using our optimality characterisations, that a CQR algorithmic variant has the optimal-order evaluation complexity of $\mathcal{O}(\epsilon^{-3/2})$ when applied to minimizing our quartically-regularised cubic subproblem, which can be further improved in special cases.  We propose practical CQR variants that judiciously use  local tensor information to construct the local cubic-quartic models. We test these variants numerically and observe them to be competitive with ARC and other subproblem solvers on typical instances and even superior on ill-conditioned subproblems with special structure.

\end{abstract}

\section{Introduction and Problem Set-up}

In this paper, we consider the unconstrained nonconvex optimization problem, 
\begin{equation}
  \min_{x\in \R^n} f(x).  
  \label{min f}
\end{equation}
Here, $f: \R^n \rightarrow \R$ is $p$-times continuously differentiable and bounded below, $p \ge 1$. 
Recent research~\cite{birgin2017worst, cartis2020sharp, cartis2020concise,  Nesterov2021implementable} showed that some optimization algorithms demonstrate superior worst-case complexity bounds when they leverage high-order derivative information of the objective function alongside adaptive regularization techniques.
In these optimization methods, a polynomial local model $m_p(x_k, s)$ plays a crucial role as an approximation to $f(x_k+s)$ at the current iterate $x_k$. Then,  $x_k$ is iteratively updated by $
    s_k \approx \argmin_{s\in \R^n}{m_p(x_k, s)}, \text{ and } x_{k+1} := x_k+s_k, 
$
whenever sufficient objective decrease is obtained. 
This process continues iteratively until an approximate local minimizer of $f$ is found, satisfying conditions such as $\|\nabla_x f(x_k)\|\leq \epsilon_g$ and $\lambda_{\min}\left(\nabla^2_x f(x_k)\right)\geq -\epsilon_H$. 

The construction of the model $m_p$ relies on the $p$th-order Taylor expansion of $f(x_k+s)$ at $x_k$
$$
   T_p(x_k, s) := f(x_k) + \sum_{j=1}^p \frac{1}{j!} \nabla_x^j f(x_k)[s]^j,
$$
where $\nabla^j f(x_k) \in \R^{n^j}$ is a $j$th-order tensor and $\nabla^j f(x_k)[s]^j$ is the $j$th-order derivative of $f$ at $x_k$ along $s \in \R^n$. 
To ensure a local model that is bounded from below and to maintain convergence in the optimization method, we introduce an adaptive regularization term\footnote{Unless otherwise stated, $\|\cdot\|$ denotes the Euclidean norm in this paper.} into $T_p$,
\begin{equation}
\tag{AR$p$ Model}
    m_p(x_k, s) = T_p(x_k, s)+ \frac{\sigma_k}{p+1}\|s\|^{p+1},
    \label{subprob}
\end{equation}
where $\sigma_k > 0$. 
Notably, when $p=1$, this equation represents the steepest descent model, while $p=2$ corresponds to a Newton-like model. In this paper, our primary focus lies on the case where $p =3$. This construction, as articulated in \eqref{subprob}, corresponds to the subproblem in the well-known adaptive regularization algorithmic framework AR$p$~\cite{birgin2017worst, cartis2020sharp, cartis2020concise}. Within AR$p$, the parameter $\sigma_k$ dynamically adjusts to ensure progress towards optimality across iterations.
Under Lipschitz continuity assumptions on $\nabla^p f(x)$, AR$p$ requires no more than $\mathcal{O}\left({\max\left[\epsilon_g^{-\frac{p+1}{p}},\epsilon_H^{-\frac{p+1}{p-1}} \right]}\right)$ evaluations of $f$ and its derivatives to compute a local minimizer to accuracy $(\epsilon_g, \epsilon_H)$ for  first-order and second-order criticality. 
 The theoretical result highlights that as we increase the order $p$, the evaluation complexity bound improves. For instance, the AR$3$ algorithm has evaluation complexity\footnote{This complexity metric excludes the computational cost associated with minimizing the subproblem \eqref{subprob}.} of $\mathcal{O}\left(\max{\epsilon_g^{-4/3}, \epsilon_H^{-2}}\right)$ which is better than the worst-case performance compared to first/second-order methods.

\subsection{Literature Review for Minimizing AR$3$ Model}

When $p = 3$, the subproblem $m_3$ is a potentially nonconvex, quartically-regularized multivariate polynomial, 
\begin{equation}
\tag{AR$3$ Model}
 m_3(x_k, s) = f(x_k) + \nabla_x f(x_k)^Ts + \frac{1}{2} \nabla_x^2 f(x_k) [s]^2 + \frac{1}{6}  \nabla_x^3 f(x_k) [s]^3+ \frac{ \sigma_k}{4} \|s\|^4.
  \label{ar3 model}
\end{equation} 
Efficiently minimizing $m_3(s)$ is important for the overall performance of the AR$3$ algorithm. Consequently, our primary objective is to devise specialized algorithms tailored for the optimization of $m_3$. While it is worth noting that finding the global minimum of $m_3$ is an NP-hard problem~\cite{burachik2021steklov, luo2010semidefinite}, our goal here is to efficiently identify a local minimizer of $m_3(s)$. Importantly, it has been established in previous work~\cite{cartis2020concise} that obtaining a local minimizer of $m_3(s)$ is sufficient to ensure the favorable computational complexity of the AR$3$ algorithm.

Finding ways to efficiently minimize $m_3$ remains an open question and is the main focus of this paper. 
While there have been some previous algorithms related to minimizing the AR$3$ subproblems, they often lack explicit handling of the tensor term or the fourth-order regularization.
Schnabel et al. \cite{chow1989derivative, schnabel1971tensor, schnabel1991tensor} considered solving unconstrained optimization using third-order tensor local models without the regularization term. 
Birgin et al. \cite{birgin2017use} introduced a variant of the AR$3$ algorithm, comparing its performance with AR$2$/ARC. 
In the case of a convex $m_3$, Nesterov proposed a series of second-order methods for its minimization \cite{Nesterov2021implementable, Nesterov2020inexact, Nesterov2021inexact, Nesterov2021superfast, Nesterov2022quartic, Nesterov2006cubic}. These methods utilize various convex optimization tools, including Bregman gradient methods, high-order proximal-point operators, and iterative minimization of convex quadratic models with quartic regularization. 
Recently, Cartis and Zhu \cite{cartis2023second, zhu2022quartic} introduced the Quadratic Quartic Regularisation (QQR) method. {QQR approximates the third-order tensor term of the AR$3$ model with a linear combination of quadratic and quartic terms, resulting in (possibly nonconvex) local models that can be solved to global optimality. In this paper, we aim to enhance the previous approach by explicitly incorporating third-order terms.} For local minimization of $m_3$, local first and second-order methods originally designed for optimizing $f(x)$ can also be applied directly, including ARC  (also known as AR$p$ with $p=2$).

In addition to the works mentioned above, there is a body of literature on polynomial optimization for (globally) minimizing certain quartic polynomials. In specific cases where the polynomial exhibits convexity or can be expressed as a sum of squares (SOS) of polynomials, semidefinite programming (SDP) methods have been shown to converge to global minimizes \cite{ahmadi2023higher, lasserre2001global, lasserre2015introduction, laurent2009sums}. Recently, \cite{ahmadi2023higher} proposed a high-order  Newton's method that uses SDP to construct and minimize an SOS convex approximation to the $p$th-order Taylor expansion of the objective function. However, when such conditions do not hold, alternative methods become necessary. In \cite{qi2004global}, a global descent algorithm is proposed for finding a global minimizer of such quartic polynomials in two variables ($n=2$).  Burachik et al. \cite{arikan2020steklov, burachik2021steklov} proposed a trajectory-type method for globally optimizing multivariate quartic normal polynomials. Another strategy involves utilizing branch-and-bound algorithms, which partition the feasible region recursively and construct nonconvex quadratic or cubic lower bounds. Despite the advancements in these methods, they all face challenges associated with the curse of dimensionality. The computational complexity can grow exponentially with the number of variables of the polynomial, rendering certain methods impractical for large-scale problems.

Therefore, our work aims to develop more efficient algorithms specifically tailored for minimizing the $m_3$ model, taking into account its tensor structure and fourth-order regularization.

\subsection{Motivation for Cubic Quartic Regularisation Method (CQR)}
\label{CQM} 

For notational simplicity\footnote{{$g^Ts$ represents the usual Euclidean inner product of $g$ and $s$, and $H[s]^2 = s^T Hs$. $T[s_1][s_2]\dotsc[s_p]$ represents the tensor $T$ applied sequentially to the vectors $s_1, s_2,\dotsc,s_p$, and similarly, $T[s]^p$ is the tensor $T$ applied repeatedly to the vector $s$ a total of $p$ times.}}, we fix $x_k$ and write $m_3(x_k, s)$ as
\begin{equation}
 m_3(s) ={{f_0 + g^Ts+ \frac{1}{2}  H [s]^2} + \frac{1}{6}  T [s]^3}+ \frac{\sigma}{4}  \|s\|^4,
 \label{m3}
\end{equation}
where $f_0 = f(x_k) =  m_3(x_k, 0) \in \R$, $g = \nabla_x f(x_k) \in \R^n$, $H = \nabla^2_x f(x_k) \in \R^{ n \times n}$ and $T = \nabla^3_x f(x_k) \in \R^{ n \times n \times n}$. 
We also denote the fourth-order Taylor expansion\footnote{For all $s^{(i)}$,  the Hessian $H_i$ is a symmetric matrix. $T_i \in \R^{n\times n \times n}$ is a supersymmetric third-order tensor which means that the entries are invariant under any permutation of its indices. {$T[s_1][s_2][s_3] =T[s_{\pi(1)}][s_{\pi(2)}][s_{\pi(3)}] $ where $\pi: \{1,2, 3\}\rightarrow\{1,2, 3\}$ is the permutation function.}}  of $m_3 (s^{(i)}+s)$ at $s^{(i)}$ as
\begin{eqnarray}
M(s^{(i)}, s)  := f_i + g_i^T s + \frac{1}{2}H_i [s]^2 +\frac{1}{6}T_i[s]^3 + \frac{\sigma}{4}\|s\|^4,
\label{taylor M}
\end{eqnarray} 
where $f_i =  m_3(s^{(i)}) \in \R$, $g_i = \nabla m_3(s^{(i)}) \in \R^n$, $H_i = \nabla^2 m_3(s^{(i)}) \in \R^{ n \times n}$ and $T_i = \nabla^3 m_3(s^{(i)}) \in \R^{ n \times n \times n}$ and $\nabla$ denotes the derivative with respect to $s$. The fourth-order Taylor expansion is exact since $m_3(s)$ is a fourth-degree multivariate polynomial. Therefore, we have $
M(s^{(i)}, s) =m_3(s+s^{(i)})$,   $ \min_{s\in \R^n}m_3(s) = \min_{s\in \R^n}M(s^{(i)}, s)$,   \text{ and } $ s^*: = \argmin_{s\in \R^n}m_3(s) = s^{(i)}+\argmin_{s\in \R^n}M(s^{(i)}, s). $

While globally minimizing the quartically regularized polynomial $M(s^{(i)}, s)$ is known to be NP-hard~\cite{luo2010semidefinite}, we employ a quadratic model with a cubic degree term and quartic regularization to provide an approximation for $M(s^{(i)}, s)$. This model, which we refer to as the \textbf{Cubic Quartic Regularization (i.e., CQR) polynomial/model} can be expressed as follows 
\begin{equation}
 M_c(s^{(i)}, s) = f_i + g_i^T s+\frac{1}{2} H_i [s]^2  +  \frac{\beta_i}{6} \|s\|_W^3+ \frac{\sigma_c^i}{4} \|s\|_W ^4 \approx M(s^{(i)}, s),
 \label{CQR Model}
\end{equation}
where $\sigma_c^i >0$, $\beta^i \in \R$, $W$ is a symmetric, positive-definite matrix\footnote{The matrix norm is defined as $\|v\|_W:= \sqrt{v^TWv}$.}. In the CQR polynomial, $\beta_i$ typically provides insights into the tensor term $T_i[s]^3$, while $\sigma_c$ is adjusted to control regularization {and algorithmic progress}. It is worth noting that our model accommodates negative values of $\beta_i$ which is particularly useful as it yields information about the 'negative tensor directions' where $T_i[s]^3 < 0$. {When $\beta_i = 0$ and $W = I_n$, $M_c(s^{(i)}, s)$ simplifies to the quadratic quartically-regularized (QQR) polynomials used in \cite{cartis2023second}. On the other hand, when $\beta_i > 0$ and $\sigma^i_c = 0$, $M_c(s^{(i)}, s)$ reduces to the quadratic model with cubic regularization polynomials, which is used in the ARC/AR$2$ algorithm~\cite{cartis2011adaptive, dussault2018arcq, kohler2017sub, martinez2017cubic, Nesterov2006cubic}. However, note that for this paper, we require $\sigma_c^i > 0$ but allow $\beta_i\leq 0$.}

The reason for choosing the CQR polynomial  $ M_c(s^{(i)}, s)$ to approximate  $M(s^{(i)}, s)$ is that in this paper, we characterse the global minimizers of such CQR polynomials and devise a procedure to locate an approximate global minimizer of $ M_c(s^{(i)}, s)$. Then, we minimize the CQR polynomial  $ M_c(s^{(i)}, s)$ iteratively and find a sequence of $\{s^{(i)}\}_{i \ge 0}$ such that $s^{(i)}$ converges to a second-order local minimizer of  $M(s^{(i)}, s)$  (or equivalently $m_3$).  We refer to this iterative minimization method as \textbf{the CQR algorithmic framework/method}. The updates for the CQR method are charactersed by the following iterations
\begin{equation*}
 s_c^{(i)} = \argmin_{s\in \R^n}  M_c(s^{(i)}, s), \qquad  s^{(i+1)} = s^{(i)} + s_c^{(i)} , \qquad i: = i+1
\end{equation*}
where $\beta_i, \sigma_c^i$ are adjusted adaptively in each iteration.

\textbf{The paper presents the following key contributions:}

\begin{itemize}
    \item We develop necessary and sufficient optimality conditions for the global minimizers of $M_c(s^{(i)}, s)$. These conditions are general and applicable to both convex and nonconvex $M_c$ and for both $\beta_i \ge 0$ and $\beta_i \le 0$. Leveraging these optimality conditions, we can transform the problem of locating a global minimizer of $M_c(s^{(i)}, s)$ into a problem of solving a univariate nonlinear equation coupled with a matrix system that is similar to the trust region and ARC secular equations (Theorem 8.2.8~\cite{cartis2022evaluation}). This opens up opportunities for developing efficient algorithms for minimizing $M_c(s^{(i)}, s)$. In this paper, we present an algorithm that is based on Cholesky factorization and Newton's root finding.

In the literature, there are related works on cubic polynomial multivariate models. For instance, Martinez and Raydan \cite{martinez2015separable, martinez2017cubic} investigate the problem of the separable cubic model, where the cubic power term is expressed as $\sum_{i=1}^n \beta_j s_j^3$, with $s_j$ being the $j$th entry of $s \in \R^n$ and $\beta_j$ being scalar constants.  {However, these studies are conducted within the framework of cubic regularization or trust region methods \cite{cartis2011adaptive, cartis2011adaptiveII,  gould2012updating, grapiglia2015convergence}. They do not incorporate third-order terms of the Taylor series in the objective function or its cubic approximations. }The minimization of quartic regularized polynomials with specific cubic terms represents a relatively new research area. As far as we are aware, there are no algorithms explicitly designed for minimizing a quadratic model with a cubic degree term and a quartic regularization term.

 \item  We then show global convergence properties of CQR algorithmic variants, analyzing the quality of the approximation provided by $M_c(s^{(i)}, s)$ to $m_3(s)$. Our theoretical analysis shows that CQR exhibits at least the same complexity as ARC and, in specific cases, the CQR algorithm demonstrates improved convergence behaviour. Practical CQR variants are proposed, where particular attention is given to the tensor approximation terms. {Preliminary numerical experiments with these CQR variants also indicate that they usually identify an approximate stationary point of $m_3$ with fewer iterations or evaluations compared to the ARC method.} 
\end{itemize}

The paper is structured as follows: In Section \ref{sec Global Optimality result for CQR}, we explore necessary and sufficient optimality conditions that hold at a global minimizer of $M_c(s^{(i)}, s)$ and present an efficient algorithm in minimizing the CQR polynomial. Section \ref{sec cqr complexity} contains the convergence proof and complexity analysis of the CQR method. In addition, Appendix \ref{sec: bound on iterations} proves that the iterates of the CQR method are uniformly bounded, a technical theorem crucial for analyzing the convergence and complexity of the CQR method. Finally, we provide algorithmic implementations and numerical examples in Section \ref{sec cqr numerical Algorithm}. The conclusion can be found in Section \ref{sec Conclusion}.

\section{Characterising the Global Optimality of the CQR Polynomial}
\label{sec Global Optimality result for CQR}

In this section, we develop necessary and sufficient conditions that hold at the global minimizers of $M_c(s^{(i)}, s)$ in \ref{CQR Model}. We convert the problem of finding a global minimizer of $M_c(s^{(i)}, s)$ into a problem of solving a nonlinear eigenvalue problem, which we can then solve using root-finding techniques. 
For the sake of notational simplicity, in this section, we write $ M_c(s^{(i)}, s)$  as $M_c(s) $, 
\begin{equation}
\label{cqr model}
 M_c(s) = f_i + g_i^T s+\frac{1}{2} H_i [s]^2  +  \frac{\beta}{6} \|s\|_W^3+ \frac{\sigma_c}{4} \|s\|_W ^4,
\end{equation}
where $\sigma_c $ represents $\sigma_c^i $ and $\beta$ as $\beta_i$.
The derivatives of $M_c(s) $ has the expression 
\begin{eqnarray}
\nabla  M_c(s)  &=& g_i + H_i s+ \frac{\beta}{2}\|s\|_W (Ws) +\sigma_c \|s\|_W^2 (Ws),
\label{gradm}
\\\nabla^2  M_c( s)  &=&  H_i + \frac{\beta}{2}\bigg(W\|s\|_W + \frac{(Ws)(Ws)^T}{\|s\|_W} \bigg) +\sigma_c \bigg(W\|s\|_W^2 + 2(Ws)(Ws)^T\bigg)
\label{hessm}
\end{eqnarray}
where $\nabla$ is the derivative with respect to $s$ and the matrix $(Ws)(Ws)^T$ is a rank-one matrix. 

The section is organized as follows: In Section \ref{sec necessary optimality}, we present  {global} necessary optimality conditions for the CQR Polynomial in \eqref{cqr model}. In Section \ref{sec Sufficient optimality}, we provide {global} sufficient optimality conditions. Finally, employing these conditions, in Section \ref{sec Newton Method}, we introduce an algorithm for locating {an approximate minimizer} of the CQR polynomial using nonlinear root-finding techniques.

\subsection{Necessary Optimality Conditions for the CQR Polynomial}
\label{sec necessary optimality}

Theorem \ref{thm: cubic opt global necessary} provides a necessary condition for the global minimizers of $M_c(s)$, and we prove this theorem using a methodology similar to that of Theorem 8.2.7 in \cite{cartis2022evaluation}.

\begin{theorem} \textbf{(Necessary optimality conditions)}
Let $s_c$ be a global minimizer of $M_c(s)$ in \eqref{cqr model} over $\R^n$ and let 
\begin{eqnarray}
    B(s_c) : = H_i + \frac{\beta}{2} W \|s_c\|_W + \sigma_c W \|s_c\|_W^2.
    \label{Bs_c}
\end{eqnarray}
Then, $s_c$  satisfies the system of equations
\begin{equation}
     B(s_c) s_c =  -g_i
    \label{nececond1}
\end{equation}
and
\begin{equation}
    B(s_c)  \succeq 0. 
    \label{nececond2}
\end{equation}
 If  $ B(s_c)$ is positive definite, $s_c$ is unique. 
\label{thm: cubic opt global necessary}
\end{theorem}

\begin{proof}
The first-order necessary local optimality conditions at $s_c$ give 
$
g_i+ H_i s_c + \frac{\beta}{2}\|s_c\|_W (Ws_c) +\sigma_c \|s_c\|_W^2 (Ws_c) = 0$ and hence that \eqref{nececond1} holds.
The second-order necessary optimality gives
\begin{equation}
w^T \bigg[ H_i + \frac{\beta}{2}\bigg(W\|s_c\|_W + \frac{(Ws_c)(Ws_c)^T}{\|s_c\|_W} \bigg) +\sigma_c \bigg(W\|s_c\|_W^2+ 2(Ws_c)(Ws_c)^T\bigg) \bigg] w \ge0 
\label{2ndopt}
\end{equation}
for all vectors $w \in \R^n$. 
If $s_c = 0$, \eqref{2ndopt} is equivalent to $ H_i$ being positive semi-definite, which immediately gives the result in \eqref{nececond2}. Thus, we only need to consider  $s_c \neq 0$.
There are two cases to consider. Firstly, suppose that $w^T W s_c = 0$. In this case, it immediately follows from \eqref{2ndopt} that
\begin{equation}
w^T \bigg[\underbrace{ H_i + \frac{\beta}{2}\|s_c\|_W W +\sigma_c \|s_c\|_W^2 W }_{=B(s_c)}\bigg] w \ge0 \text{ for all } w \text{ for which }w^T W s_c  = 0. 
\label{result1}
\end{equation}
It thus remains to consider vectors $w$ for which $w^T W s_c \neq 0$. Since $w$ and $s_c$ are not orthogonal with respect to the $W$-norm, the line $s_c + \tilde{k} w$ intersects the ball of radius $\|s_c\|$ at two points, $s_c$ and $u_*$ and thus
\begin{equation}
    \|s_c\|_W = \|u_*\|_W.
    \label{eqnorm}
\end{equation}
We set $w_* = u_*- s_c$, and note that $w_*$ is parallel to $w$.
Since $s_c$ is a global minimizer, we immediately have that
\begin{eqnarray*}
    0 &\le& M_c(u_*) - M_c(s_c) 
    \\&=& g_i^T(u_* - s_c) +  \frac{1}{2} {u_*}^T  H_i u_* - \frac{1}{2} {s_c}^T H_i s_c+\frac{\beta}{6} (\|u_*\|_W^3-\|s_c\|_W^3) +\frac{\sigma_c}{4}(\|u_*\|_W^4-\|s_c\|_W^4)
    \\&=& g_i^T(u_* - s_c) +  \frac{1}{2} {u_*}^T  H_i u_* - \frac{1}{2} {s_c}^T H_i s_c
\end{eqnarray*}
where the last equality follows from \eqref{eqnorm}. But \eqref{nececond1} gives 
\begin{eqnarray}
(u_* - s_c)^T  g_i= ( s_c - u_*)^T H_i s_c + \frac{\beta}{2} \|s_c\|_W( s_c - u_*)^TW s_c+ \sigma_c \|s_c\|_W^2 ( s_c - u_*)^TW s_c.
 \label{gus}
\end{eqnarray}
In addition, \eqref{eqnorm} gives 
\begin{eqnarray}
( s_c - u_*)^TWs_c = \frac{1}{2}{s_c}^TWs_c + \frac{1}{2}{u_*}^TWu_* - {u_*}^TWs_c = \frac{1}{2}{w_*}^TWw_*=  \frac{1}{2} \|w_*\|_W^2.
 \label{eqnorm2}
\end{eqnarray}
Thus we have
\begin{eqnarray*}
    0 &\le& M_c(u_*) - M_c(s_c) 
    \\ &\overset{\eqref{gus}}{=}& ( s_c - u_*)^T H_i s_c + \frac{\beta}{2} \|s_c\|_W( s_c - u_*)^T Ws_c+ \sigma_c \|s_c\|_W^2 ( s_c - u_*)^T Ws_c  +  \frac{1}{2} {u_*}^T  H_i u_* - \frac{1}{2} {s_c}^T H_i s_c
  \\  &\overset{\eqref{eqnorm2}}{=}& ( s_c - u_*)^T H_i s_c + \frac{\beta}{4}  \|s_c\|_W\|w_*\|_W^2+ \frac{\sigma_c}{2} \|s_c\|_W^2 \|w_*\|_W^2  +  \frac{1}{2} {u_*}^T  H_i u_* - \frac{1}{2} {s_c}^T H_i s_c
  \\  &=& \frac{1}{2}{ w_*}^T  H_i w_* + \frac{\beta}{4}   \|s_c\|_W\|w_*\|_W^2+ \frac{\sigma_c}{2} \|s_c\|_W^2 \|w_*\|_W^2 =  \frac{1}{2} {w_*}^T\bigg[ H_i +  \frac{\beta}{2}  \|s_c\|W+\sigma_c \|s_c\|_W^2W \bigg]{w_*}. 
\end{eqnarray*}
We deduce that
\begin{equation}
w^T \bigg[ H_i + \frac{\beta}{2}\|s_c\|_W W +\sigma_c \|s_c\|_W^2 W \bigg] w \ge0 \text{ for all } w \text{ for which }w^T W s_c  \neq 0. 
\label{result2}
\end{equation}
Hence, \eqref{result1} and \eqref{result2} together show that $H_i + \frac{\beta}{2}\|s_c\|_W W +\sigma_c \|s_c\|_W^2 W \succeq 0$. The uniqueness of $s_c$ when $H_i + \frac{\beta}{2}\|s_c\|_W W +\sigma_c \|s_c\|_W^2 W$ is positive definite follows immediately from \eqref{nececond1}.

\end{proof}

\subsection{Sufficient Optimality Conditions for the CQR Polynomial}
\label{sec Sufficient optimality}
In this section, we derive sufficient optimality conditions using two approaches. Specifically, in Section \ref{subsec: sufficient Nocedal}, we establish a sufficient optimality condition inspired by Nocedal and Wright \cite[Thm 4.1]{nocedal1999numerical}, while in Section \ref{sec: secular eqn}, inspired by Cartis et al \cite[Sec 8.2.1]{cartis2022evaluation}, we deduce another set of sufficient optimality conditions. Each approach corresponds to different cases of sufficiency. Finally, in Section \ref{subsec: integrate two appro}, we integrate both proofs to provide an overall sufficient theorem (Theorem \ref{thm case by case sufficiency}) that can be applied to all CQR polynomials.

\subsubsection{Sufficient Optimality Conditions using Quadratic Forms Reformulations}
\label{subsec: sufficient Nocedal}
In this section, we investigate sufficient conditions that hold at a global minimizer of  $M_c(s)$, following a methodology inspired by Nocedal and Wright \cite[Thm 4.1]{nocedal1999numerical}.

\begin{theorem}  \textbf{(Sufficient optimality conditions using quadratic forms reformulations)} 
Let $ B(s_c)$ be defined as in \eqref{Bs_c}. Then, $s_c$ is the global minimizer of $M_c(s)$ in \eqref{cqr model} over $\R^n$ if the following conditions are satisfied:
\begin{enumerate}
\item $g$ is in the range of $ B(s_c)$, such that
\begin{equation}
    B(s_c) s_c = \bigg(H_i + \frac{\beta}{2} W \|s_c\|_W I_n + \sigma_c W \|s_c\|_W^2 I_n\bigg) s_c = -g_i.  
    \label{Suffcond1}
\end{equation}
\item $ B(s_c)$ is positive semidefinite, such that
\begin{equation}
   B(s_c) := H_i + \frac{\beta}{2} W \|s_c\|_W + \sigma_c W \|s_c\|_W^2 \succeq 0. 
    \label{Suffcond2}
\end{equation}
\item $\beta$ satisfies
\begin{eqnarray}
   \beta \ge  -3\sigma_c\|s_c\|_W ,\qquad  \text{or equivalently} \qquad  \|s_c\|_W \ge -\frac{ \beta}{3  \sigma_c}. 
    \label{condition on beta}
\end{eqnarray}
\end{enumerate}
The global minimizer is unique if $B(s_c)$ is positive definite or $\beta >  -3\sigma_c\|s_c\|_W$. 
\label{thm: cubic sufficient Nocedal}
\end{theorem}

\begin{proof}
We assume wlog $f_i $ = 0. For all vectors $w \in \R^n$, we have
\small{}
\begin{eqnarray*}
    M_c(s_c+w) & =& g_i^T (s_c+w)  +\frac{1}{2} H_i [s_c+w]^2   +    \frac{\beta}{6} \|s_c+w\|_W^3+ \frac{\sigma_c}{4} \|s_c+w\|_W ^4,
    \\  & =&  \underbrace{g_i^T (s_c+w) + \frac{1}{2}B(s_c) [s_c+w]^2 - \frac{\beta}{12} \|s_c\|_W^3 - \frac{\sigma_c}{4} \|s_c\|_W^4}  _{\mathcal{F}_1} 
     \\&&  \underbrace{-\frac{1}{2} \bigg[ \frac{\beta}{2}  \|s_c\|_W + \sigma_c  \|s_c\|_W^2 \bigg] \|s_c+w\|_W^2  +  { \frac{\beta}{6} \|s_c+w\|_W^3+ \frac{\sigma_c}{4} \|s_c+w\|_W ^4  } +  {\frac{\beta}{12} \|s_c\|_W^3 + \frac{\sigma_c}{4} \|s_c\|_W^4}}_{\mathcal{F}_2}. 
\end{eqnarray*}
\normalsize{}
After rearranging, 
\begin{eqnarray*}
\mathcal{F}_1  = \bigg[\underbrace{g_i^T s_c + \frac{1}{2} B(s_c) [s_c]^2 {- \frac{\beta}{12} \|s_c\|_W^3 - \frac{\sigma_c}{4} \|s_c\|_W^2}}_{= M_c(s_c)}\bigg] + \bigg[ \underbrace{g_i  + B(s_c)s_c}_{= 0 \text{ by } \eqref{Suffcond1}} \bigg]^T  w+ \underbrace{\frac{1}{2} B(s_c)[w]^2}_{ \ge 0 \text{ by } \eqref{Suffcond2}}. 
\end{eqnarray*}
Therefore, we deduce that $\mathcal{F}_1 \ge M_c(s_c)$. Also, if $B(s_c)$ is positive definite, then,  $\mathcal{F}_1 > M_c(s_c)$ for all $w \in \R^n$ with equality only at $w =0$. 

On the other hand, we denote $E =  \|s_c\|_W$ and $F =  \|s_c+w\|_W$, clearly $E \ge 0$ and  $F\ge 0$. We can simplify and rearrange $\mathcal{F}_2$ as
\begin{eqnarray}
\mathcal{F}_2 
= \frac{1}{2}(E-F)^2 \bigg[\frac{\beta}{6}(E+2F)  + \frac{\sigma_c}{2} (E+F)^2 \bigg]. 
\label{f2 ineq}
\end{eqnarray}
Clearly, when $\beta \ge 0$, we have $\mathcal{F}_2 \ge 0$.  {Else, when $\beta \le 0$, we let $\tilde{f}(F) := \frac{(E+F)^2}{(E+2F)}$ for $F \ge 0$; then $ F_*: =\argmin_{F \ge 0} \tilde{f}(F)  = 0$ and $\tilde{f}(F_*)  =E$. This implies that $E$ is the best lower bound for $\tilde{f}(F)$ with $F \ge 0$. Together with condition  \eqref{condition on beta},} we have
$
  0 \ge \beta \ge -3 \sigma_c E \ge - 3\sigma_c \frac{(E+F)^2}{(E+2F)}. 
$
This immediately gives $\frac{\beta}{6}(E+2F)  + \frac{\sigma_c}{2} (E+F)^2 \ge 0$ and, consequently, $\mathcal{F}_2 \ge 0$. Note that if $\beta >  -3 \sigma_c\|s_c\|_W$,   equality is attained only at $w=0$. 

Combining  $\mathcal{F}_1 \ge 0$ and  $\mathcal{F}_2 \ge 0$, we have $M_c(s_c+w)  \ge M_c(s_c) $ for all $w \in \R^n$ and $s_c$ is a global minimizer of $M_c$. If $B(s_c)$ is positive definite or $\beta >  -3 \sigma_c\|s_c\|_W$, we have $M_c(s_c+w) > M_c(s_c) $ for all $w \in \R^n/\{0\}$. 

\end{proof}

\begin{remark} \textbf{(Discrepancy between sufficient and necessary conditions)}
The sufficiency conditions \eqref{Suffcond1}--\eqref{Suffcond2} coincide with the necessary conditions \eqref{nececond1}--\eqref{nececond2}. The distinction between sufficiency and necessity arises from the inclusion of \eqref{condition on beta}, which imposes certain constraints on the choice of $\beta$. {Namely, $\beta$ should not be selected to have an excessively negative value. } 
\end{remark}

\begin{remark} \textbf{(Reduction to a cubically regularized quadratic polynomial)}
Theorem \ref{thm: cubic sufficient Nocedal} can be considered a generalization of sufficient conditions for the minimization of a cubically regularized quadratic polynomial \cite{cartis2011adaptive, Nesterov2006cubic}. 
If $\sigma_c=0$ and $\beta = 2 \sigma_q > 0$, we recover the sufficient condition for the cubically regularized quadratic polynomial,
$$
M_{q}(s)  = f_i + g_i^T s + \frac{1}{2} H_i [s]^2  +  \frac{\sigma_q}{3} \|s\|_W ^3. $$
The sufficiency condition \eqref{condition on beta} is naturally satisfied since $\beta > 0$. The sufficiency conditions \eqref{Suffcond1} and \eqref{Suffcond2} for the CQR polynomial reduce to the sufficiency conditions of the cubic regularized polynomial in  \cite[Thm 8.2.8]{cartis2022evaluation}. 
\end{remark}

\subsubsection{Sufficient Optimality Conditions using a Secular Equation Approach}
\label{sec: secular eqn}

\begin{theorem} \textbf{(Sufficient optimality conditions using a secular equation approach)} Assume that $ H_i$ and $\beta$ satisfy any one of the following cases, 
\begin{itemize}
    \item \textbf{Case 1: }$ H_i$ is indefinite.  
    \item \textbf{Case 2: }$ H_i$ is positive definite and $\beta \ge 0$.
    \item \textbf{Case 3: }$ H_i$ is positive definite, $ -4 \sqrt{\sigma_c\lambda_1} < \beta \le 0$, 
    and $\big\|\big[H_i -\frac{\beta^2}{16\sigma_c} W\big]^{-1}g\big\|_W\le -\frac{\beta}{4\sigma_c}$ where $\lambda_1$ is the leftmost eigenvalue of the pencil $( H_i ; W)$. 
\end{itemize}
Let $s_c \in \R^n$ satisfy \eqref{nececond1}--\eqref{nececond2}. Then, $s_c$ is a global minimizer of $M_c(s)$ in \eqref{cqr model} over $\R^n$.
\label{thm:sufficient secular equation}
\end{theorem}

Theorem \ref{thm:sufficient secular equation} shares some similarities with Theorem \ref{thm: cubic sufficient Nocedal}. In particular, when $\beta \ge 0$, both proofs indicate that conditions \eqref{nececond1}--\eqref{nececond2} represent the necessary and sufficient condition for global optimality of $s_c$ for $M_c$. When $\beta<0$, \eqref{condition on beta} in Theorem \ref{thm: cubic sufficient Nocedal} establishes the relationship between $\beta$ and $\|s_c\|_W$ required to satisfy the sufficient optimality conditions. On the other hand, Case 3 in Theorem \ref{thm:sufficient secular equation} provides the relationship between $\beta$ and the eigenvalues of $H_i$ to identify the global minimizer before computing $\|s_c\|_W$. A combined result of the two sufficiency theorems can be found in Section \ref{subsec: integrate two appro}. In this subsection, our focus is on proving Theorem \ref{thm:sufficient secular equation}.

To prove Theorem \ref{thm:sufficient secular equation}, we first need to introduce the \textbf{secular equation}. To satisfy the conditions \eqref{nececond1}--\eqref{nececond2}, we seek a vector $s(\lambda)$ which satisfies the system. 
\begin{tcolorbox}[breakable, enhanced, title = Secular Equation]
\vspace{-0.5cm}
\begin{eqnarray}
B(s_c) s(\lambda) &=& (H_i + \lambda W)s(\lambda) = -g_i , 
\label{cond0 nonlinear}
\\ B(s_c) &=& H_i + \lambda W \succeq 0, 
\label{cond1 nonlinear}
\\ \lambda  &=& \frac{\beta}{2}\|s(\lambda)\|_W+\sigma_c \|s(\lambda)\|_W^2.
\label{cond2 nonlinear}
\end{eqnarray}
where $B(s_c)$ is defined in \eqref{Bs_c}. 
\end{tcolorbox}

To rearrange \eqref{cond1 nonlinear}, let $\lambda_1 \le \lambda_2 \le \dotsc \le \lambda_n$ be the eigenvalues of the pencil $( H_i ; W)$ and let $u_1, u_2, \dotsc, u_n$ be the corresponding generalised eigenvectors\footnote{More details can be found in Appendix A.5 of \cite{cartis2022evaluation}.} where each $u_i \in \R^n$. Then,
\begin{eqnarray}
    H_i U = WUD, \text{ where } U^T W U = I_n,
    \label{general eig value 1}
\end{eqnarray}
where $D \in \R^{n \times n}$ is the diagonal matrix with entries  $\lambda_i$ for $1 \le i \le n$, and $U \in \R^{n \times n}$ is the matrix whose columns
are the $u_i$ vectors for $1 \le i \le n$. Define $g_u:=U^T g_i$ and $s_u(\lambda)=U^TWs(\lambda)$. From this, we may deduce that  $s(\lambda) = U s_u(\lambda)$ and $g = W Ug_u$ using \eqref{general eig value 1}. It then follows from \eqref{cond1 nonlinear} that 
\begin{eqnarray*}
   (H_i + \lambda W) s(\lambda) =  (H_i + \lambda W) U s_u(\lambda) \underset{\eqref{general eig value 1}}{=} WU(D + \lambda I)  s_u(\lambda)  = -g_i = -WU g_u. 
\end{eqnarray*}
Since $W$ and $U$ are nonsingular, therefore $ s_u(\lambda):= (D + \lambda I)^{-1} g_u. $ For $\lambda \neq -\lambda_i$ for any $1 \le i \le n$, 
\begin{eqnarray}
 \Psi(\lambda) := \|s(\lambda)\|_W^2 =\|s_u\|^2_2 = \|(D + \lambda I )^{-1}U^T g_i\|^2 = \sum_{i=1}^n \frac{\gamma_i^2}{(\lambda + \lambda_i)^2}
\label{general eig value}
\end{eqnarray}
where $\gamma_i$ is the $i$th component of $U^T g_i$. 

We also rearrange \eqref{cond2 nonlinear} to a nonlinear function of $\lambda$, such that
\begin{equation}
   k_\pm (\lambda) := \frac{1}{4\sigma_c}(-\beta \pm\sqrt{\beta^2 + 16\lambda \sigma_c}) = \|s(\lambda)\|_W. 
    \label{two trenches}
\end{equation}
We say \eqref{two trenches} is well defined if $k_\pm (\lambda) = \|s(\lambda)\|_W \ge 0$. Using  \eqref{general eig value} and \eqref{two trenches}, we give the nonlinear equations that need to be satisfied for $\lambda$. 
\begin{tcolorbox}[breakable, enhanced, title = Nonlinear Equations for $\lambda$]
Let $\psi(\lambda) = \sqrt{\Psi(\lambda)}$, where $\Psi(\lambda)$ is defined in \eqref{general eig value}.  Then, solving the secular equations \eqref{cond0 nonlinear}--\eqref{cond2 nonlinear} is equivalent to finding roots that satisfy
\begin{equation}
 \psi(\lambda)  = k_\pm (\lambda) = \frac{1}{4\sigma_c}(-\beta \pm \sqrt{\beta^2 + 16\lambda \sigma_c}) \quad \text{for} \quad \lambda \ge -\lambda_1.
    \label{secular} 
\end{equation}
where $\lambda_1$ be the leftmost eigenvalue of the pencil $( H_i; W)$.
\end{tcolorbox}

The sufficiency implication states that if we can prove the existence of a solution $\lambda^*$ to \eqref{secular}, then it follows that the solution to the system \eqref{cond0 nonlinear}--\eqref{cond2 nonlinear} exists.
{
\begin{itemize}
\item If $H_i + \lambda^* W \succ 0$, \eqref{cond0 nonlinear} yields a unique $s_c$, indicating a unique global minimizer for $M_c(s)$.
\item The hard case, similar to the trust-region subproblem \cite[Sec 8.3.1]{cartis2022evaluation}, may occur when $\tau_1 := u_1^T g_i = 0$, leading to non-unique global minimizers for $M_c(s)$.
\end{itemize}
}

\textbf{Case-by-case proof: global minimizer of $M_c(s)$ for $\tau_1 = u_1^T g_i\neq 0$.}

We start with the standard cases when $g_i$ is not orthogonal to the {first} generalised eigenvector of the pencil $( H_i; W)$ (i.e., $\tau_1 \neq 0$). {We  prove that, in these cases, there exist $\lambda^*$ satisfying \eqref{secular}, and $ H_i + \lambda^* W $ is positive definite. Thus, $M_c(s)$ has a unique global minimizer. To aid the proof, we define point $\mathcal{A}$ as the intersection point for $k_+$ and $k_-$ which has the coordinate  $(\lambda_A, s_A): = \big(-\frac{\beta^2}{16\sigma_c}, -\frac{\beta}{4\sigma_c}\big)$.} 

\begin{figure}[ht]
    \centering
    \includegraphics[width = 15.5cm]{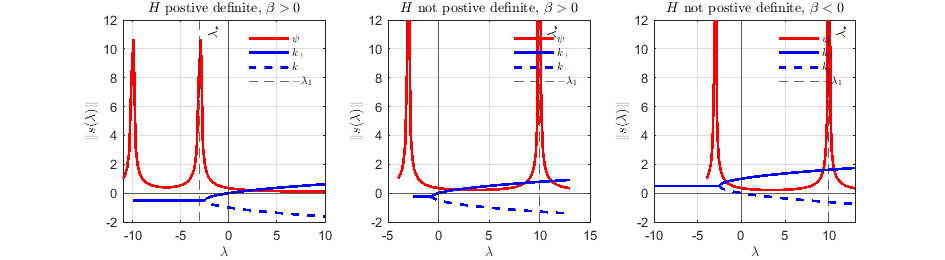}
  \caption{\small  $\sigma_c = 10$, $W = I_2$, $g_i = [1,-1]^T$. \textbf{Left plot (Case 2):} $\beta = 20$, $H_i  = [10, 0; 0, 3]$; 
  \textbf{Mid plot (Case 1):} $\beta = 20$, $H_i = [-10, 0; 0, 3]$;  
  \textbf{Right plot (Case 1):} $\beta = -20$, $H_i  = [-10, 0; 0, 3]$. }
    \label{fig beta>0}
\end{figure}

\textbf{Case 1:} If $ H_i$ is indefinite and $\tau_1 \neq 0$, we seek a solution satisfying $\lambda > -\lambda_1 \ge 0$. As illustrated by the right two plots in Figure \ref{fig beta>0}, for all $\beta \in \R$ and $\lambda > -\lambda_1 \ge  0$, only $k_+(\lambda)$ is well defined (i.e.  $k_+(\lambda)  \ge  0$).  Since $\psi(\lambda)  =\sqrt{\Psi(\lambda)}$ is monotonically decreasing with  $\lim_{\lambda \rightarrow -\lambda_1}\psi(\lambda) \rightarrow  \infty$ and $\lim_{\lambda \rightarrow \infty }\psi(\lambda) \rightarrow  0$ and  $k_+(\lambda)$ is monotonically increasing with $k_+(-\lambda_1) \ge 0$, there is a unique solution $\lambda^*$ in the $k_+(\lambda)$ trench of \eqref{secular}. {Note that $\lambda^* > -\lambda_1$ and $H_i + \lambda^* W \succ 0$, ensuring the existence of a unique global minimizer.}

\textbf{Case 2:} {If $\beta \ge 0$, $H \succ 0$ and $\tau_1 \neq 0$, the intersection point $\mathcal{A} := (\lambda_A, s_A) = \big(-\frac{\beta^2}{16\sigma_c}, -\frac{\beta}{4\sigma_c}\big)$ must be in the third quadrant and the only trench that is well defined is the $k_+(\lambda)$ for $\lambda > 0$ (An example illustrating this is given in the left plot in Figure \ref{fig beta>0})}. Moreover, $k_+(\lambda)$ is monotonically increasing for $\lambda > 0$ with $k_+(0) = 0$ and $\lim_{\lambda \rightarrow \infty }k_+(\lambda) \rightarrow  \infty $. $\psi(\lambda)$ is monotonically decreasing with $\psi(0) > 0$ and $\lim_{\lambda \rightarrow \infty }\psi(\lambda) \rightarrow  0$. Thus, there is a unique solution $\lambda^*$  lies in the $k_+(\lambda)$ trench for \eqref{secular} {with  $\lambda^* > -\lambda_1$. The uniqueness of the global minimizer follows similarly to Case 1.}

\begin{figure}[ht]
    \centering
    \includegraphics[width = 15.5cm]{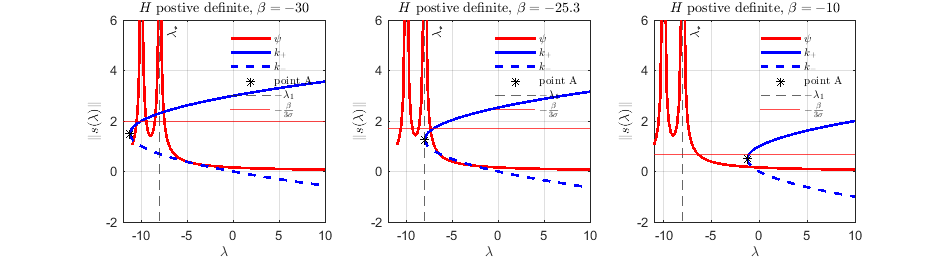}
 \caption{\small  $\sigma_c = 5$, $W = I_2$, $g_i  = [1,-1]^T$, $H_i  = [10, 0; 0, 8]$. \textbf{Left plot  (Case 4):} $\beta = -30$; 
  \textbf{Mid plot (Case 5):} $\beta = -25.3 $;  
  \textbf{Right plot (Case 3):} $\beta = -10$. } 
    \label{fig beta<0}
\end{figure}

\textbf{Case 3:} If $\beta \leq 0$, $H \succ 0$, and $\tau_1 \neq 0$, as shown in Figure \ref{fig beta<0}, \eqref{secular} may have multiple solutions, so we need to impose extra conditions on $\beta$. 
 The unique solution is in the $k_-$ trench if point $\mathcal{A}$ stays at the right of $\lambda = -\lambda_1$ and point $\mathcal{A}$ is above $\psi$ (the third plot of Figure \ref{fig beta>0}). In other words, 
$$\lambda_A = -\frac{\beta^2}{16\sigma_c} > -\lambda_1, \qquad \text{and} \qquad s_A = -\frac{\beta}{4\sigma_c} >\psi\bigg(-\frac{\beta^2}{16\sigma_c}\bigg) . $$
{Note that we have $\lambda^* \ge \lambda_A > -\lambda_1$, therefore $H_i + \lambda^* W \succ 0$, ensuring a unique global minimizer.}

\begin{lemma} \textbf{(Case 3)} Assume $H \succ 0$, and $\tau_1 \neq 0$,  if $\beta \le 0$ and satisfies $\frac{\beta^2}{16\sigma_c} < \lambda_1$ and $\psi(-\frac{\beta^2}{16\sigma_c}) = \big\|\big[H_i -\frac{\beta^2}{16\sigma_c} W\big]^{-1}g\big\|_W < -\frac{\beta}{4\sigma_c}$, then there exists a unique solution $\lambda^*$ to \eqref{secular} that lies in the $k_-(\lambda)$ trench. 
\label{lemma case 3}
\end{lemma}

\begin{proof}

Since $\beta \le 0$, we have $\lambda_A = -\frac{\beta^2}{16\sigma_c} <0$ and $s_A = \frac{\beta}{4\sigma_c}>0$. Point $\mathcal{A}$ which has coordinate   $(\lambda_A, s_A): = (-\frac{\beta^2}{16\sigma_c}, -\frac{\beta}{4\sigma_c})$ is in the second quadrant.  
Since $k_-$ is only well-defined for $\lambda > \lambda_A$, the range of interest for finding a solution is $\lambda \in (\lambda_A, 0)$. Both $\psi$ and $k_-$ are continuous and monotonic in the interval $(\lambda_A, 0)$, and $\psi(\lambda_A) \le -\frac{\beta}{4\sigma_c} = k_-(\lambda_A)$ while $\psi(0) > 0 = k_-(0)$, while $\psi(0) > 0 = k_-(0)$, there must be a unique intersection between $k_-$ and $\psi$ within this interval. 
Lastly, we can clarify that the feasible root cannot occur in the $k_+$ trench. This is because $k_+$ is monotonically increasing and $\psi$ is monotonically decreasing for all $\lambda > \lambda_A$ and $k_+(\lambda_A) = k_-(\lambda_A) > \psi(\lambda)$. 
\end{proof}

\textbf{Case-by-case proof: non-unique global minimizer of $M_c(s)$ for $\tau_1 = u_1^T g_i=0$.}


Now we consider the case where $g_i$ is orthogonal to the generalised eigenvector of the pencil $( H_i; W)$ (i.e., $\tau_1 = u_1^T g_i=0$). 
If $\tau_1 = 0$, we may have difficulties similar to the hard case of the trust-region subproblem \cite[Sec 8.3.1]{cartis2022evaluation}. 

Under these circumstances, \eqref{cond0 nonlinear} will be consistent if we set $\lambda = \lambda_s := -\lambda_1$, such that \eqref{cond0 nonlinear}  becomes $
     (H  +\lambda_s W) s_s = -g_i. $ 
{Thus, $H+\lambda_s W \succeq 0$, so \eqref{cond1 nonlinear} is satisfied. }However, since $\lambda_1 = -\lambda_s$ is the leftmost eigenvalue of the pencil $( H_i ; W)$, it follows that $(H_i + \lambda_s W) u_1 = 0$ for all $u_1$ corresponding to $\lambda_1$ and thus
\begin{eqnarray}
    g_i + (H_i +\lambda_s W) (s_s+k u_1) =0
\end{eqnarray}
for any scale $k \in \R$. Thus, in order to satisfy \eqref{cond2 nonlinear}, we must additionally choose the scalar $k = k_s$ such that $
    \lambda_s  = \frac{\beta}{2}\big\|s_s+k_su_1\big\|_W+\sigma_c\big\|s_s+k_su_1\big\|_W^2.$ 
To obtain an exact solution in the hard case requires the generalised eigenvalue $\lambda_1$, a corresponding eigenvector $u_1$ of $ H_i$, and the `trajectory' vector $
s_s = \lim_{\lambda \rightarrow \lambda_s} s(\lambda ) = -(H_i + \lambda_s W)^+g$ 
where $(H_i + \lambda_s W)^+$ is the generalised inverse of $H_i + \lambda_s W$. We summarize these results in the Theorem \ref{thm hard case}. More details and full proof can be found in \cite[Corollary 8.3.1]{cartis2022evaluation}.

\begin{theorem} 
Any global minimizer of $M_c$ may be expressed as
$$
    s_c = s_s\text{ uniquely}, \text{ if }   \tau_1 \neq 0, \qquad \text{ or}
\qquad s_c=  s_s + k_s u_1,  \text{ if }  \tau_1 = 0,
$$
where $ \tau_1 = u_1^T g_i$, $
    s_s = -(H_i +\lambda_{s} W)^+ g_i$ and $\lambda_s  = \frac{\beta}{2}\|s_c\|_W+\sigma_c \|s_c\|_W^2 > -\lambda_1$, where $k_s$ is one of the two roots of $
    \lambda_s  = \frac{\beta}{2}\big\|s_s+k_su_1\big\|_W+\sigma_c\big\|s_s+k_su_1\big\|_W^2. $ 
\label{thm hard case}
\end{theorem}


\subsubsection{General Sufficient Optimality Conditions for the Global Minimum of $M_c$ }
\label{subsec: integrate two appro}

In this subsection, we integrate Theorem \ref{thm: cubic sufficient Nocedal} and Theorem \ref{thm:sufficient secular equation} to provide a general sufficient optimality result along with a graphical representation of the feasible values of $\beta$ (Figure \ref{fig region plot}).
The sufficient conditions in Theorem \ref{thm:sufficient secular equation} correspond to Cases 1 to 3 in Theorem \ref{thm case by case sufficiency}. The sufficiency conditions in Theorem \ref{thm: cubic sufficient Nocedal} not only confirm the results of Theorem \ref{thm:sufficient secular equation} but also extend the range of cases that we can address (see Case 4).

\begin{tcolorbox}[breakable, enhanced, title = General Sufficient Optimality Conditions for the Global Minimum of $M_c$] 
\begin{theorem}
{Consider the CQR polynomial $M_c(s)$ in \eqref{cqr model}, let $s_c \in \R^n$ be a vector that satisfies \eqref{cond0 nonlinear}, \eqref{cond1 nonlinear}, and \eqref{secular} with the corresponding trench selected as described in Case 1 to Case 5 below. Then, $s_c$ is a global minimizer of $M_c(s)$ over $\R^n$. }
\begin{itemize}
    \item \textbf{Case 1:} (From Theorem \ref{thm: cubic sufficient Nocedal} and Theorem \ref{thm:sufficient secular equation}) $H_i$ is indefinite, take $k_+$ trench in \eqref{secular}. 
    \item \textbf{Case 2:} (From Theorem \ref{thm: cubic sufficient Nocedal} and Theorem \ref{thm:sufficient secular equation}) $H_i$ is positive definite and $\beta \ge 0$, take $k_+$ trench in \eqref{secular}. 
    \item \textbf{Case 3:} (From Theorem \ref{thm:sufficient secular equation}) $ H_i$ is positive definite,
    \begin{eqnarray}
   -4 \sqrt{\sigma_c\lambda_1} \le \beta \le 0, \quad \text{and} \quad  -\frac{\beta}{4\sigma_c} >\psi\bigg(-\frac{\beta^2}{16\sigma_c}\bigg)  =\bigg\|\bigg(H_i -\frac{\beta^2}{16\sigma_c} W\bigg)^{-1}g\bigg\|_W, 
   \label{k- condition}
   \end{eqnarray} 
   take $k_-$ trench in \eqref{secular}.
    \item \textbf{Case 4:} (From Theorem \ref{thm: cubic sufficient Nocedal}) $H_i$ is positive definite,
    $
\beta \le  -3\sqrt{2}\sqrt{\lambda_1 \sigma_c} \approx -4.24 \sqrt{\lambda_1 \sigma_c}, 
$
take $k_+$ trench in \eqref{secular}.
\item  \textbf{Case 5: }In the rest of the cases, namely, $H_i$ is positive definite and $ -3\sqrt{2}\sqrt{\lambda_1 \sigma_c} \le \beta \le -4\sqrt{\lambda_1 \sigma_c}$. Or $H_i$ is positive definite, 
\begin{eqnarray*}
-4 \sqrt{\lambda_1 \sigma_c} \le \beta \le 0   \quad \text{and}    \quad    -\frac{\beta}{4\sigma_c} <\psi\big(-\frac{\beta^2}{16\sigma_c}\big)  =\bigg\|\big(H_i -\frac{\beta^2}{16\sigma_c} W\big)^{-1}g\bigg\|_W. 
\end{eqnarray*}
Up to three solutions satisfying \eqref{secular} may exist. The solution that minimizes the function value of $M_c$ corresponds to the global minimizer. A root, denoted by $(\lambda_+, s_+)$, within the $k_+$ trench exists, and this solution satisfies $M_c(s_+)<M_c(0)$.
\end{itemize}
\label{thm case by case sufficiency}
\end{theorem}
\end{tcolorbox}

\begin{remark} 
Cases 1, 2, and 3 follow from Theorem \ref{thm:sufficient secular equation}. Cases 1, 2, and 4 are derived from Theorem \ref{thm: cubic sufficient Nocedal}, and therefore, $\|s_c\|_W$ satisfies $\|s_c\|_W \ge -\frac{\beta}{3\sigma_c}$.
\end{remark}

\begin{remark} 
{All values of $\beta \in \R$ are addressed in Theorem \ref{thm case by case sufficiency}. This indicates that by finding the root of \eqref{secular} that satisfies \eqref{cond0 nonlinear} and \eqref{cond2 nonlinear}, we can find a global minimizer of $M_c$ for all $\beta \in \R$ and $\sigma_c > 0$.
The intersection point of the $k_+$ and $k_-$ trenches, denoted as $\mathcal{A} = (\lambda_A, s_A) = \left(-\frac{\beta^2}{16\sigma_c}, -\frac{\beta}{4\sigma_c}\right)$, corresponds to various cases outlined in Theorem \ref{thm case by case sufficiency}. When $H$ is indefinite, regardless of the location of point $\mathcal{A}$, we are in Case 1 (as shown in the left plot of Figure \ref{fig region plot}). If $H$ is indefinite, when point $\mathcal{A}$ falls within specific colored regions, it corresponds to the cases described in Theorem \ref{thm case by case sufficiency} (as shown in the right plot of Figure \ref{fig region plot}).}
As illustrated in Figure \ref{fig region plot}, Cases 1--4 encompass the majority of $\beta$ values, while Case 5 covers the remaining situations. An example of Case 5 is presented in the second plot of Figure \ref{fig beta<0}.
 
\begin{figure}[ht]
    \centering
    \includegraphics[width = 15.5cm]{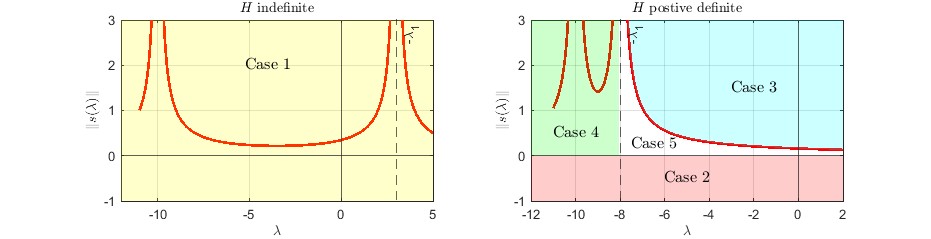}
    \caption{\small Illustration plot depicting the location of point $\mathcal{A}$ corresponding to different cases outlined in Theorem \ref{thm case by case sufficiency}.
     $\sigma_c = 5$, $W = I_2$, $g_i  = [1,-1]^T$, $H_i  = [10, 0; 0, -3]$ for the left plot and  $H_i  = [10, 0; 0, 8]$ for the right plot. The red curve represents $\psi(\lambda)$. 
   }
    \label{fig region plot}
\end{figure}
\end{remark}

\begin{proof}
{\textbf{Case 1 and Case 2: follows from Theorem \ref{thm: cubic sufficient Nocedal} and Theorem \ref{thm:sufficient secular equation}}.} According to Theorem \ref{thm: cubic sufficient Nocedal}, if $\beta \ge 0$, we can observe from \eqref{f2 ineq} that $\mathcal{F}_2 > 0$, and consequently, $M_c(s+w) > M_c(s)$ for all $w \in \R^n$. The condition on $\beta$ is naturally satisfied.
According to Theorem \ref{thm:sufficient secular equation}, as discussed in the case-by-case analysis, there exists a unique root $\lambda^* > -\lambda_1$ in the $k_+(\lambda)$ for the nonlinear equation \eqref{secular}.

{\textbf{Case 3: follows from Theorem \ref{thm:sufficient secular equation}}.} This case is established using Theorem \ref{thm:sufficient secular equation} and Lemma \ref{lemma case 3}. It is worth noting that this case cannot be derived using Theorem \ref{thm: cubic sufficient Nocedal}. The reason is that any solution found in the $k_-$ trench satisfies $\|s_c\|_W \le s_A=-\frac{\beta}{4\sigma_c}$, which does not always meet the conditions of \eqref{condition on beta}.

{\textbf{Case 4: follows from Theorem \ref{thm: cubic sufficient Nocedal}.}} Assuming $\beta \leq 0$, $H \succ 0$, and $\tau_1 \neq 0$, the sufficiency conditions provided in Theorem \ref{thm: cubic sufficient Nocedal} offer an alternative approach to characterse the global minimizer. If point $\mathcal{A}$ is positioned to the left of $\lambda = -\lambda_1$ and the value of $k_+$ at $\lambda = -\lambda$ satisfy $k_+(-\lambda_1) \geq -\frac{\beta}{3\sigma_c}$ (as depicted in the first plot of Figure \ref{fig beta<0}), then there exists a solution that satisfies the sufficient condition in Theorem \ref{thm: cubic sufficient Nocedal}. This solution lies within the $k_+$ trench and represents the global minimizer of $M_c$. In other words, to meet the sufficient condition in Theorem \ref{thm: cubic sufficient Nocedal}, the following conditions are required:
\begin{eqnarray*}
\lambda_A = -\frac{\beta^2}{16\sigma_c} \le -\lambda_1, \quad \text{and} \quad k_+(-\lambda_1) =\frac{1}{4\sigma_c}(-\beta + \sqrt{\beta^2 -16\lambda_1 \sigma_c}) \ge -\frac{\beta}{3\sigma_c}.
\end{eqnarray*}
This implies $\beta \le -3\sqrt{2}\sqrt{\lambda_1 \sigma_c}$.

{\textbf{Case 5:}} In the remaining cases, there may be up to three solutions satisfying \eqref{secular} as depicted in the second plot of Figure \ref{fig beta<0}. By the necessary optimality conditions, the solution that minimizes the function value of $M_c$ corresponds to the global minimizer. As point $\mathcal{A}$ either lies below the curve $\psi(\lambda)$ or to the left of $\lambda = -\lambda_1$, a root in the $k_+$ trench is always guaranteed. 
Let $(\lambda_+, s_+)$  be the root in the $k_+$ trench satisfying \eqref{secular} and consequently satisfying \eqref{cond0 nonlinear}--\eqref{cond2 nonlinear}. We deduce that
\begin{eqnarray*} 
 M_c(0)-    M_c( s_+) &=& - g_i^T s_+ -  \frac{1}{2} H_i [s_+]^2  -  \frac{\beta}{6} \|s_+\|^3 - \frac{\sigma_c}{4}\|s_+\|^4  \underset{\eqref{cond0 nonlinear}}{=}  \frac{1}{2} H_i [s_+]^2 +  \frac{\beta}{3} \|s_+\|^3+ \frac{3\sigma_c}{4}\|s_+\|^4   
 \\&\ge&   \|s_+\|^2 \bigg(\frac{1}{2} \lambda_1 +  \frac{\beta}{3} \|s_+\|+ \frac{3\sigma_c}{4}\|s_+\|^2 \bigg) >0.
\end{eqnarray*}
Since $\lambda_1>0$ and $-3\sqrt{2}\sqrt{\lambda_1 \sigma_c} \le \beta < 0$, the discriminant of the quadratic equation $q(\|s_+\|):= \frac{1}{2} \lambda_1 +  \frac{\beta}{3} \|s_+\|+ \frac{3\sigma_c}{4}\|s_+\|^2 $ is negative, i.e.,  $\frac{\beta^2}{9} - \frac{3\lambda_1\sigma_c}{2} < 0.$ This gives $ M_c(0)-    M_c( s_+)  \ge 0$. 
\end{proof}

\subsection{An Algorithm for Minimizing the CQR Polynomial}
\label{sec Newton Method}

{As suggested by Theorem \ref{thm: cubic opt global necessary} and Theorem \ref{thm case by case sufficiency}, we have transformed the task of identifying a global minimizer of $M_c$ into the problem of solving the nonlinear equations \eqref{secular} for $\lambda$ in the corresponding trench.
To numerically determine the root of \eqref{secular}, we utilize the global optimality charactersations outlined in Sections \ref{sec necessary optimality} to \ref{sec Sufficient optimality}. Additionally, we extend well-established Cholesky factorization-based root-finding techniques, as detailed in \cite[Ch 8.2.3]{cartis2022evaluation}.
The primary objective of this subsection is to formulate an algorithm for locating the root of the univariate nonlinear equation. We opt to employ univariate Newton iteration to find the necessary root of \eqref{secular}. Further convergence analysis for the root-finding method is deferred to future work. }

In practice, as suggested by \cite[Sec 6.1]{cartis2011adaptive}, instead of solving \eqref{secular}, it may be preferable to solve 
\begin{equation}
  \phi_1(\lambda) := \big[\psi(\lambda)\big]^{-1}-  \big[k_{\pm}(\lambda)\big]^{-1} = 0 \qquad \text{for} \qquad \lambda \ge  \lambda_1
  \label{phi_1}
\end{equation}
where the choice between using $k_+$ trench or $k_-$ trench is made on a case-by-case basis, as determined by Theorem \ref{thm case by case sufficiency}. 
The updates for the Newton iteration are detailed in Lemma \ref{Newton lemma}.


\begin{lemma}
Suppose $g \neq 0$, 
the Newton iteration updates are 
\begin{eqnarray}
    \Delta \lambda^{k} =  \frac{\phi_1(\lambda)} {\phi_1'(\lambda)} = \frac{\|s(\lambda)\|_W^{-1} - \big[k_{\pm}(\lambda)\big]^{-1}}{\|\omega(\lambda)\|_W^2\|s(\lambda)\|_W^{-3} - \big[(k_{\pm}(\lambda))^{-1}\big]'},
    \qquad \lambda^{(k+1)} = \lambda^{(k)}+     \Delta \lambda^{k}
     \label{Newton correction}
\end{eqnarray}
where $\omega (\lambda)$ is given by $W^{-1}(H_i + \lambda W) = L(\lambda) L^T(\lambda)$ and $L(\lambda) \omega (\lambda) = s(\lambda)$. Also, 
\begin{eqnarray*}
k_{\pm}^{-1}(\lambda) = \frac{2 \sigma_c }{-\beta  \pm \sqrt{\beta^2 + 4 \lambda \sigma_c}}, \quad  \big[(k_{\pm}(\lambda))^{-1}\big]' = \frac{\mp 4 \sigma_c^2}{ (-\beta \pm \sqrt{\beta^2+4\lambda\sigma_c})^{2}(\beta^2+4\lambda\sigma_c)^{1/2} }.
\end{eqnarray*}
\label{Newton lemma}
\end{lemma}

\noindent
\textit{Sketch of Proof.}
Based on \cite[Lemma 6.1]{cartis2011adaptive}, we have given the Cholesky factorization $W^{-1}(H_i + \lambda W) = L(\lambda) L^T(\lambda)$, then
\begin{eqnarray}
    \Psi'(\lambda) = 2 s(\lambda)^T W \nabla_{\lambda}  s(\lambda) =  -2\|\omega(\lambda)\|_W^2
\end{eqnarray}
where $L(\lambda) L^T(\lambda) s(\lambda) = -Wg$ and $L(\lambda) \omega (\lambda) = s(\lambda)$. Consequently, 
\begin{eqnarray}(\psi(\lambda)^{-1})' =(\Psi(\lambda)^{-1/2})' = - \frac{1}{2}\Psi(\lambda)^{-3/2} \Psi'(\lambda)  =  \|s(\lambda)\|_W^{-3}\|\omega(\lambda)\|_W^2.
\label{derivative}
\end{eqnarray}
Substituting \eqref{derivative} into the Newton updates $\Delta \lambda^{k} =  \frac{\phi_1(\lambda)} {\phi_1'(\lambda)}$ gives the result. More details are in Appendix \ref{derivative appendix} and Lemma 7.3.1 in \cite{conn2000trust}.

\textbf{Numerical Set-up and Initialization for the Newton Step}

{In a similar manner to the techniques used to solve the trust region and ARC subproblems \cite[Ch 8.2.3]{cartis2022evaluation} \cite[Sec 6.1]{cartis2011adaptive},} for the minimization of the CQR polynomial in \eqref{CQR Model}, we employ Newton's method to solve the secular equation \eqref{phi_1}. As described in Theorem \ref{thm case by case sufficiency}, the solution of the secular equation \eqref{phi_1} falls into either Case 1, 2, 4, or 5, where it resides in the $k_+$ trench, or Case 3, where it lies within the $k_-$ trench.

In practical implementation, we first search for a solution in the $k_+$ trench using the initialization in Remark \ref{remark initialize Newton}. If no solution is found in the $k_+$ trench, we proceed to search for a solution in the $k_-$ trench. Based on our numerical experiments in Section \ref{sec cqr numerical Algorithm}, we find that in approximately 90\% of the iterations, the solution lies in the $k_+$ trench. 
The algorithm is detailed in Algorithm \ref{Newton algo}, and the initialization details are provided in Remark \ref{remark initialize Newton}

\begin{algorithm}
\caption{\small Newton's method to Solve \eqref{phi_1}}
\textbf{Initialization}: Initialize the counter $k=0$, $s^{(0)}=0$  and $\lambda^{(0)}$ as in Remark \ref{remark initialize Newton}. 
 Set an accuracy level $\epsilon_N>0$, usually $\epsilon_N = 10^{-8}$ smaller than the accuracy of the subproblem. 
\\\textbf{Initialize trench:} Choose $k_+$ trench: $k_T(\cdot) := k_+ (\cdot)$. 
\\\textbf{Input:} $f_i$,  $g_i$, $H_i$, $\beta$ 
 and $\sigma_c>0$. 
 
\While{$\left|\|s^{(k)}\| - k_T(\lambda^{(k)})\right| >\epsilon$}
{\eIf {$H_i + \lambda^{(k)} I_n \succ 0$}
{Cholesky factorize $(H_i + \lambda^{(k)}  I_n) = LL^T$. 
Solve $LL^Ts^{(k)}=- g_i$ and $L\omega^{(k)} = s^{(k)}$. 
\\Compute the Newton update step:
\begin{eqnarray*}
    \Delta \lambda^{k} =  \frac{\phi_1(\lambda^{(k)})} {\phi_1'(\lambda^{(k)})} = \frac{\|s^{(k)}\|_W^{-1} - [k_T(\lambda^{(k)})]^{-1}}{\|\omega^{(k)}\|_W^2\|s^{(k)}\|_W^{-3} - \big[(k_T(\lambda^{(k)}))^{-1}\big]'},
    \qquad \lambda^{(k+1)} = \lambda^{(k)}+     \Delta \lambda^{k}.
\end{eqnarray*}
}
{Switch to $k_-$ trench: $k_T(\cdot) = k_-$. Set $\lambda^{(0)} = 0$. }
}
\label{Newton algo}
\end{algorithm}

\begin{tcolorbox}[breakable, enhanced, title = Initialization for the Newton Algorithm]
\begin{remark}
In Cases 1, 2, 4, and 5, where the solution resides in the $k_+$ trench, we initialize $\lambda^{(0)}$ as follows: 
\begin{itemize}
    \item Case 1 and Case 2: If $H_i$ is indefinite, we set $\lambda^{(0)} := \max\{-\lambda_1, 0\} + 10^{-6}$. This ensures $H + \lambda^{(0)} I \succ 0.$
    \item Case 4 and 5:  If $H_i  \succ 0$ and $\beta < 0$, we set $\lambda^{(0)} := \max \{ -\lambda_1, \lambda_A\} + 10^{-6}$ where $\lambda_A = -\frac{\beta^2}{16\sigma_c}$ denotes the intersection point for $k_+$ and $k_-$.  This initialization also ensures $H + \lambda^{(0)} I \succ 0.$ This choice is motivated by the observation that, as illustrated in the two leftmost plots of Figure \ref{fig beta<0}, the solution to the secular equation \eqref{phi_1} consistently resides in the right side of $\lambda = -\lambda_{\min}[H]$ and $\lambda = \lambda_A$.
\end{itemize}
If no solution is found in the $k_+$ trench, it implies we are in Case 3. Consequently, we proceed to seek a solution within the $k_-$ trench. In Case 3, where $H_i \succ 0$ and $\beta < 0$, as illustrated on the rightmost plot in Figure \ref{fig beta<0}, the solution to the secular equation \eqref{phi_1} is near $\lambda = 0$ and falls within the interval $[\lambda_A, 0]$. Therefore, we initialize $\lambda^{(0)} = 0$. 
\label{remark initialize Newton}
\end{remark}
\end{tcolorbox}

\begin{remark}
Regardless of the case, we ensure that $H_i + \lambda^{(0)} I_n$ is positive definite, making the Cholesky factorization possible in the first iteration. While it is worth noting that this does not entirely avoid the hard case, as it is only the first iteration, this is a measure we take to ensure that we do not start the algorithm with an indefinite $H_i + \lambda^{(0)} I_n$. Throughout the iterations of Algorithm \ref{Newton algo}, in all the numerical examples we examined, the positive definiteness of $H_i + \lambda^{(k)} I_n$ is preserved, and we did not encounter the hard case. The Cholesky factorization also remains valid in the cases we tested. However, note that Algorithm \ref{Newton algo} is a preliminary implementation of Newton root-finding for the nonlinear equation \eqref{secular}. Addressing and implementing the hard case, as well as the convergence analysis of Algorithm \ref{Newton algo}, is dedicated to future work.
\end{remark}

\begin{remark}
    
{Theorem \ref{thm: cubic opt global necessary} and Theorem \ref{thm case by case sufficiency} give a global optimality 
characterization for the CQR polynomial in \eqref{CQR Model}. 
However, in Case 5 of Theorem \ref{thm case by case sufficiency},  up to three solutions may exist. Then,
since we are using Newton's method, depending upon the initialization, Algorithm \ref{Newton algo} may converge to a root in $k_+$ which corresponds to a local minimizer of $M_c$. To locate the global minimizer in Case 5, one needs to locate all the roots in the $k_+$ and $k_-$ trench (there are at most three roots) and evaluate the value of $M_c$ at these roots. The roots that give the lowest function value of $M_c$ will correspond to  global minimizers of $M_c$. We defer the implementation of this additional root-finding procedure to future algorithmic developments.} 
Numerically, in all the cases that we tested, the algorithm converges to a global minimizer in Cases 1--4 and converges to a local minimizer $s^{(K)}$, with a decrease in $M_c$ value in Case 5, such that $M_c(s^{(K)}) < M_c(s^{(0)})$. 
\end{remark}


Our current approach here, for solving the system \eqref{cond0 nonlinear}--\eqref{cond2 nonlinear} and  \eqref{secular}, uses Cholesky factorization and Newton's method to address the secular equation \eqref{secular}, and is just one of many possibilities.
Similarly to existing approaches for trust region and regularization subproblems (\cite[Ch. 6]{cartis2011adaptive} and \cite[Ch. 8, 10]{cartis2021scalable}), scalable iterative algorithms could be designed here, based on Krylov methods, eigenvalue formulations, or subspace optimization. While we recognize the potential of these tools for solving \eqref{cond0 nonlinear}--\eqref{cond2 nonlinear} and finding the global minimum of the CQR polynomial, our primary focus in this paper is on using the CQR polynomial to minimize the $m_3$ subproblem. As a result, we defer a detailed analysis of these alternative approaches to future work.

\section{The CQR Algorithmic Framework for Minimizing $m_3$}
\label{sec cqr complexity}

The CQR algorithm minimizes the AR$3$ subproblem (i.e., $m_3$ as defined in \eqref{m3}) by globally minimizing a sequence of potentially nonconvex CQR polynomials (i.e., quadratic models with cubic term and quartic regularization){\footnote{For a fixed subproblem in \eqref{m3}, $x_k$ is fixed, and we have access to the coefficients for the subproblem; namely, $f_0 = f(x_k)$, $g = \nabla_x f(x_k)$, $H = \nabla^2_x f(x_k)$, $T = \nabla^3_x f(x_k)$, and $\sigma>0$.}.} {Additionally,  the CQR algorithmic framework is tensor-free in the sense that we do not require the entries of $T$ to be stored; we only need vector and matrix products with the tensor $T$ (such as $T[s]$ or $T[s]^2$).
We define the norms of the derivatives as $\|g\| := \sqrt{g^Tg}$, $\|H\| := \max_{i}|\lambda_{i}(H)|$, and $\|T\|=\max_{\|u_1\|=\|u_2\|=\|u_3\|=1} \big| T [u_1][u_2][u_3] \big| < \infty$.}
In this section and the rest of the paper, when constructing the CQR model $M_c$ in \eqref{CQR Model}, we choose $W = I_n$. Algorithm \ref{TCQR} outlines the framework of the CQR algorithm. 

\begin{algorithm}
\caption{\small CQR Algorithmic Framework with Adaptive Regularization}
\textit{Initialization}: Set $s^{(0)} = \boldsymbol{0} \in \R^n$ and $i=0.$ An initial regularization parameter $d_0 = 1$, constants $ \eta_1 > \eta > 0$, $\gamma > 1 > \gamma_2 >0$.
\\
\textit{Input}: $f_0:=m_3(s^{(0)})$,  $g_0 := \nabla m_3(s^{(0)})$, $H_0 : = \nabla^2 m_3(s^{(0)})$, $\sigma>0$ is the given regularization parameter in $m_3$ as defined in \eqref{m3};  an accuracy level $\epsilon>0$. 

\textbf{Step 1: Test for termination.} If $\|g_i\| = \|\nabla m_3(s^{(i)})\| \le \epsilon$, terminate with the $\epsilon$-approximate first-order minimizer $s_\epsilon = s^{(i)}$. 

\textbf{Step 2: Step computation.}
Compute $\beta_i$ to give a scalar representation of the (local) tensor information in $m_3$\footnotesize{$^a$}. \normalsize{}Set
\begin{equation}
\tag{CQR Model}
M_c(s^{(i)}, s) : = f_i + g_i s + \frac{1}{2}  H_i [s]^2 +  \frac{1}{6} \beta_i \|s\|^3+\bigg(\underbrace{\frac{\sigma}{4}+d_i}_{:=\sigma_c^i}\bigg) \|s\|^4.
\label{t mc}
\end{equation}
Compute $s_c^{(i)} = \argmin_{s \in \R^n} M_c(s^{(i)}, s)$ (see Remark \ref{remark solver s_c}). 

\textbf{Step 3: Acceptance of trial point.} Compute $m_3(s^{(i)}+s_c^{(i)})$ and define
\begin{equation}
\tag{Ratio Test}
\rho_i = \frac{m_3(s^{(i)}) -m_3(s^{(i)}+s_c^{(i)})}{m_3(s^{(i)}) - M_c(s^{(i)}, s_c^{(i)})}.
\label{ratio test}
\end{equation}
\textbf{Step 4: Regularization Parameter Update.} 
\begin{itemize}
\item If $\rho_i \ge \eta_1$ \textit{(very successful step)}: Set
\\$s^{(i+1)} = s^{(i)} + s_c^{(i)} $, $d_{i+1} = \gamma_2 d_i$, update $\beta_{i+1}$ (details in Section \ref{sec cqr numerical Algorithm}), and compute $g_{i+1}$. 
\item If $ \eta \le \rho_i \le \eta_1$ \textit{(successful step)}: Set
\\$s^{(i+1)} = s^{(i)} + s_c^{(i)} $, $d_{i+1} =  d_i$,  update $\beta_{i+1}$, and compute $g_{i+1}$.
\item  If $\rho_i \le \eta$, \textit{(unsuccessful step)}\\$s^{(i+1)} = s^{(i)}$, increase $d_{i+1} = \gamma d_i$, $\beta_{i+1} = \beta_i$.
\end{itemize}
Increment $i$ by one and go to Step 1. 
\\
\footnotesize{$^a$Details about selecting $\beta_i$ can be found in Algorithm \ref{PCQR} in Section \ref{sec cqr numerical Algorithm}.}
\label{TCQR}
\end{algorithm}
\renewcommand\footnoterule{}

\begin{remark}
    \textbf{(Purpose of $\beta_i$)} The purpose of $\beta_i$ is to provide a scalar approximation of the third-order information $\frac{T[s]^3}{\|s\|^3}$. For a given tensor $T \in \R^{n \times n \times n}$, $\frac{T[s]^3}{\|s\|^3}$ is bounded. Also, the term $\frac{T[s]^3}{\|s\|^3}$ can be either positive or negative. Therefore, we allow $\beta_i$ to take both positive and negative values. {We also assume that $\beta_i$ is bounded for all iterations, as outlined in Assumption \ref{assumption1}. 
     The bound for $\beta_i$ is used for the convergence of the CQR algorithm, ensuring optimal complexity and practical numerical implementation. 
     This bound can be easily applied in practice, such as by setting $B = \max_{1 \le \iota, j, k \le n} \big|T[ \iota, j, k]\big|$ in the practical CQR algorithm (Algorithm \ref{PCQR} in Section \ref{sec cqr numerical Algorithm}), to prevent $|\beta_i|$ from becoming excessively large.}  More details about the computation of $\beta_i$ are discussed in Section \ref{sec cqr numerical Algorithm}. 
    \label{purpose of di} 
\end{remark}

\begin{assumption} 
For every iteration, $i$, $\beta_i$ is uniformly bounded, such that $|\beta_i| \le B$.
\label{assumption1}
\end{assumption}

\begin{remark}
\label{purpose of beta}
    \textbf{(Purpose of $d_i$)}  $d_i$ is an adaptive regularization parameter; $d_i \ge 0$ ensures that $M_c$ remains bounded from below.
In some iterations, if $\beta_i$ is positive, then $\beta_i \|s\|^3$ can serve as a regularization term for the $M_c$ model, potentially superseding the role of the quartic regularization term $d_i \|s\|^4$. Thus, in these iterations, we have more flexibility in choosing $d_i$ (i.e., we can set $d_i = 0$).
\end{remark}

\begin{remark} {(\textbf{Minimization of $M_c$}) 
In Algorithm \ref{TCQR}, we assume that we find a global minimizer of $M_c$, such that $s_c^{(i)} = \argmin_{s \in \R^n} M_c(s^{(i)}, s)$. 
For optimal complexity analysis, we focus on a variant of the CQR algorithmic framework (Section \ref{sec: Convergence and Complexity of CQR}). While our proof in Section \ref{sec: Convergence and Complexity of CQR} (in particular Lemma \ref{lemma: lower bound of norms}) assumes exact minimization of $\nabla M_c(s^{(i)}, s_c) = 0$, similar bounds can be derived when considering approximate minimization of $M_c$ with a step-termination condition, such as $\|\nabla M_c(s^{(i)}, s_c)\| \le \Theta \|s_c\|^3$ for some $\Theta \in (0, 1)$. Such step-termination conditions are also used in \cite{cartis2020concise}, with alternative variants available \cite{cartis2011adaptive}.}
\label{remark solver s_c}
\end{remark}

It is worth noting that the CQR method can be applied directly to the minimization of an objective function (not just $m_3$). If the objective function satisfies suitable Lipschitz conditions and is bounded below, then the CQR method can also be used for local minimization of a potentially nonconvex objective function. However, our interest here is not to create a general method for optimizing $f$, instead, we focus on efficiently minimizing $m_3$.

This section is organized as follows: In Section \ref{sec: Convergence and Complexity of CQR}, we introduce a variant of the CQR algorithmic framework and establish its convergence guarantee. We demonstrate that this variant requires at most $O(\epsilon^{-3/2})$ function and derivative evaluations to compute an approximate first-order critical point of $m_3$ that satisfies $\|g_i\| \le \epsilon$. This bound is at least as good as that of the ARC algorithm for minimizing $m_3$. In Section \ref{sec: CQR Algorithm in Special Cases}, we explore specific instances of $m_3$ where CQR algorithm exhibits improved convergence behavior.

\subsection{Convergence and Complexity of a CQR Algorithmic Variant}
\label{sec: Convergence and Complexity of CQR}

Algorithm \ref{algo: cqr variant 1} is a variant of the CQR algorithmic framework (Algorithm \ref{TCQR}) {for which we prove an optimal worst-case complexity bound. }

\begin{algorithm}
\caption{\small A CQR Variant of Algorithm \ref{TCQR}}
Fix a constant $\alpha \in (0, \frac{1}{2})$. Initialize and proceed through Steps 1 to 3 following the same procedure as outlined in the CQR algorithmic framework (Algorithm \ref{TCQR}). In Step 4, we utilize the following rule for updating the regularization parameter, $\beta_{i+1}$ and $s^{(i+1)}$. 

\If{$\|\nabla m_3(s^{(i)}+s_c^{(i)})\| \le \epsilon$}
{{\textbf{Terminate} with an $\epsilon$-approximate first-order minimizer $s_\epsilon$ of  $m_3(s^{(i)}+s_c^{(i)})$.}}
\Else{
\uIf{$\rho_i \ge \eta$ and $\beta_i \in [\alpha, B] $} 
{\textit{{Successful Step}}: $s^{(i+1)} := s^{(i)} + s_c^{(i)} $, update $0 \le  d_{i+1} \le d_i$, update $\beta_{i+1}$ (details in Section \ref{sec cqr numerical Algorithm}), compute $g_{i+1}$. }
\uElseIf{$\rho_i \ge \eta$ and $\beta_i \in [-B, -4\alpha] $ and $\sigma_c^i  \notin  \big[\frac{1}{6}, \frac{2}{3}\big] \big(-\beta_i + \alpha\big)\|s_c^{(i)}\|^{-1}$}
{\textit{{Successful Step}}: $s^{(i+1)} := s^{(i)} + s_c^{(i)} $, update $0 \le  d_{i+1} \le d_i$, update $\beta_{i+1}$, compute $g_{i+1}$. }
\uElseIf{$\rho_i \ge \eta$ and $\beta_i \in [-4\alpha, \alpha] $ and \begin{eqnarray}
    \sigma_c^i \ge \frac{2}{3} \big(-\beta_i + \alpha\big)\|s_c^{(i)}\|^{-1}
    \label{addition cond on di}
\end{eqnarray}} 
{\textit{{Successful Step}}: $s^{(i+1)} := s^{(i)} + s_c^{(i)} $, update $0 \le  d_{i+1} \le d_i$, update $\beta_{i+1}$, compute $g_{i+1}$. }
\Else{
\textit{\textbf{Unsuccessful Step}}: $s^{(i+1)} := s^{(i)}$,  $\beta_{i+1}:= \beta_{i}$ and increase $d_{i+1} := \gamma \max\{1, d_{i}\}. $
}
}
\label{algo: cqr variant 1} 
\end{algorithm}

{In this variant, we introduce a new constant, $\alpha \in [0, \frac{1}{2}]$, to further control the step length $s_c$. Due to this additional parameter and the step size control, we need to consider two additional cases in which we reject the step and increase regularization. 
The first case occurs when $-B \le \beta_i \le -4 \alpha$, and $\sigma_c$ falls within the narrow region of $\big[\frac{1}{6}, \frac{2}{3}\big] \big(-\beta_i + \alpha\big)|s_c^{(i)}|^{-1}$. The second case arises when $\beta_i$ is in the narrow range of $[-4 \alpha, \alpha]$, and $\sigma_c^i \ge \frac{2}{3} \big(-\beta_i + \alpha\big)|s_c^{(i)}|^{-1}$.
For the remaining cases where $\beta_i$ and $\sigma_c^i$ do not fit into these two scenarios, the ratio test defined in \eqref{ratio test} with an outcome of $\rho_i \ge \eta$ is sufficient to determine the success of the step.}

It is also important to note that whenever we satisfy  $\rho_i \ge \eta$ as illustrated in the CQR algorithmic framework (Algorithm \ref{TCQR}), even without the requirement for $\sigma_c^{(i)}$ in Algorithm \ref{algo: cqr variant 1}, we are guaranteed to have a decrease in values of $m_3$ at each successful iteration, 
$$
 m_3(s^{(i)}) -m_3(s^{(i+1)}) \ge \eta \bigg[m_3(s^{(i)}) - M_c(s^{(i)}, s_c^{(i)}) \bigg] >0
$$
where  $\eta \in (0, 1]$. In essence, the CQR algorithmic framework (Algorithm \ref{TCQR}) ensures convergence to an approximate local minimum of $m_3$ that adheres to first-order optimality conditions. {But our aim in this section is to establish a global rate of complexity bound, hence the need to introduce an additional requirement involving $\alpha.$}

To establish convergence and determine the complexity bound for Algorithm \ref{algo: cqr variant 1}, we adopt a framework similar to the convergence proof for AR$1$ as presented in \cite[Sec 2.4.1]{cartis2022evaluation}. The proof is structured as follows:
\begin{enumerate}
    \item \textbf{Groundwork}: We establish upper and lower bounds for the step size (Theorem \ref{thm upper bound on step size} and Lemma \ref{lemma: lower bound of norms}) using an upper bound on the tensor term (Corollary \ref{corollary upper bound for Ti}).
    \item \textbf{Decrease in the value of $m_3$}: We prove that the CQR model constructed in successful steps of Algorithm \ref{algo: cqr variant 1} yields a decrease in the value of $m_3$ by at least $O(\epsilon^{-3/2})$, where $\epsilon$ represents the prescribed first-order optimality tolerances (Theorem \ref{thm cubic Upper Bound for s^i} and Theorem \ref{thm M3 value decrease}).
    \item \textbf{Bounded Iterations}: We establish that there is a limited number of unsuccessful iterations (Lemma \ref{lemma upper bound on di} and Lemma \ref{lemma Number of Iterations}).
\end{enumerate}
Combining these results, we derive the complexity bound for Algorithm \ref{algo: cqr variant 1} (Lemma \ref{lemma Bound on the number of successful iterations} and Theorem \ref{thm: cqr complexity}). Specifically, Algorithm \ref{algo: cqr variant 1} is guaranteed to find a first-order critical point of $m_3$ within at most $\mathcal{O}(\epsilon^{-3/2})$ function and derivatives' evaluations.

 Under Assumption \ref{assumption1}, Theorem \ref{thm upper bound on step size} and Lemma \ref{lemma: lower bound of norms} establishes a uniform upper and lower bound for the norm of the iterates $s^{(i)}$ for both successful and unsuccessful iterations.

\begin{theorem}
\textbf{(An upper bound on step size)} Suppose that  Assumption \ref{assumption1}  holds and Algorithm \ref{algo: cqr variant 1} is employed. Then, the norm of the iterates (both successful and unsuccessful) are uniformly bounded above independently of $i$, such that 
$\|s^{(i)}\|<r_c$ for all $i \ge 0$, where 
$r_c$ is a constant that depends only on the coefficients of $m_3$. 
\label{thm upper bound on step size}
\end{theorem} 
\begin{proof}
{The proof of Theorem \ref{thm upper bound on step size} uses a similar technique to \cite[Lemma 3.2]{cartis2019universal}.} The proof for Theorem \ref{thm upper bound on step size} is given in of Appendix \ref{sec: bound on iterations}.
\end{proof}

\begin{corollary}
 \textbf{(Upper bounding  the tensor term)} 
Let $L_H: = \Lambda_0 +6 \sigma r_c$, where $$\Lambda_0 :=\max_{\|u_1\|=\|u_2\|=\|u_3\|=1} {T [u_1][u_2][u_3]}.$$ Then, for every iteration, $i$, 
\begin{eqnarray*}
  \|T_i\|:=  \max_{\|u_1\|=\|u_2\|=\|u_3\|=1} \big|T_i [u_1][u_2][u_3] \big|= \Lambda_0 +6 \sigma \|s^{(i)}\| < \Lambda_0 +6 \sigma r_c = L_H.
\end{eqnarray*}
\label{corollary upper bound for Ti}
\end{corollary}
\begin{proof}
See Lemma 3.3 \cite{cartis2023second}. 
\end{proof}
Note that $L_H$ is an iteration-independent constant is a constant that depends only on the coefficients in $m_3$. 

\begin{lemma} 
\label{lemma: lower bound of norms}
\textbf{(A lower bound on the step size)} Suppose that  Assumption \ref{assumption1}  holds and Algorithm \ref{algo: cqr variant 1} is used. Assume that $\|\nabla m_3(s^{(i)} + s_c^{(i)})\| > \epsilon$ and $\|g_i\| > \epsilon$, then if \footnote{If $d_i = 0$, $\|s_c^{(i)}\| > (B  +L_H )^{-1/2}\epsilon^{1/2}$.} $d_i \ge 0$, 
\begin{eqnarray}
\|s_c^{(i)}\| >   \min\bigg\{(B  +L_H )^{-1/2}\epsilon^{1/2}, \frac{1}{2}{d_i}^{-1/3}  \epsilon^{1/3} \bigg\}.
\label{lower bound on step size}
\end{eqnarray}
\end{lemma}

\begin{proof} 
For notational simplicity, in this proof, we use $s_c$ to represent $s_c^{(i)}$. 
Since $s_c$ is the global minimum of $M_c$, we have $\nabla M_c(s^{(i)}, s_c) = 0$. Thus, 
$$
\epsilon<\|\nabla m_3(s^{(i)} + s_c)\| = \| \nabla  M(s^{(i)}, s_c)\| = \|\nabla M_c(s^{(i)}, s_c) - \nabla  M(s^{(i)}, s_c)\|. 
$$
Using the expression of $\nabla M_c$ in \eqref{gradm} and  $\nabla M(s^{(i)}, s) = g_i + H_i s +\frac{1}{2}T_i[s]^2 + \sigma\|s\|^2 s$, we have
\begin{eqnarray}
\epsilon< \bigg\| \frac{1}{2} \beta_i \|s_c\| s_c -  \frac{1}{2}  T_i [s_c]^2 +   4d_i   \|s_c\|^2 s_c \bigg\| \le \frac{B}{2}\|s_c\|^2 +    \frac{L_H}{2} \big\|s_c\big\|^2+ 4d_i\|s_c\|^3
\label{diff grad}
 \end{eqnarray}
 where the last inequality uses norm properties, $ \beta_i  > -B$ from Assumption \ref{assumption1} and $\big\|T_i\big\| \le L_H$ from Corollary \ref{corollary upper bound for Ti}. 
If $d_i \ge 0$, either one of the following must be true
\begin{eqnarray*}
\frac{ \epsilon}{2} < \bigg(\frac{B}{2}  +  \frac{L_H}{2}  \bigg) \|s_c\|^2 \Rightarrow   \|s_c\| > (B+ L_H)^{-1/2}\epsilon^{1/2},
\quad \text{or} \quad  \frac{ \epsilon}{2} <  4 d_i \|s_c\|^3   \Rightarrow\|s_c\| >\frac{1}{2}{d_i}^{-1/3}  \epsilon^{1/3}
\end{eqnarray*}
which gives \eqref{lower bound on step size}. 
\end{proof}

Lemma \ref{lemma: lower bound of norms} demonstrates that we can establish a lower bound for the norm of the iterates $s^{(i)}$.  The lower bound on step size plays a crucial role in guaranteeing a decrease in values of $m_3$. The proof incorporates standard techniques found in \cite[Lemma 2.3]{birgin2017worst} and in \cite[Lemma 3.3]{cartis2020concise}, which demonstrate that the step cannot be arbitrarily small compared to the criticality conditions. It is important to note that Lemma \ref{lemma: lower bound of norms} applies for both successful and unsuccessful iterations provided that the criticality conditions are not met, i.e.,  {$\|\nabla m_3(s^{(i)} + s_c^{(i)})\| > \epsilon$ and $\|g_i\| \ge \epsilon$. }

Theorem \ref{thm cubic Upper Bound for s^i} provides an upper bound for $m_3$ {at $m_3(s^{(i)}+s_c^{(i)}) = M(s^{(i)},s_c^{(i)})$ where $m_3$ and $M$ are defined in \eqref{m3} and \eqref{taylor M}, respectively. } 

\begin{theorem}
\label{thm cubic Upper Bound for s^i}
\textbf{(A local upper bound on $m_3$)} Suppose that  Assumption \ref{assumption1}  holds and Algorithm \ref{algo: cqr variant 1} is used.
Then, the following results are true. 
\begin{enumerate}
    \item \textbf{\bb{(Case 1)}} If $\beta_i > L_H$, for any $d_i \ge 0$, we have
    $ 
     m_3(s^{(i)}+s_c^{(i)}) = M(s^{(i)},s_c^{(i)}) \le M_c(s^{(i)},s_c^{(i)}). 
    $
    \item \textbf{\bb{(Case 2)}} If any $\beta_i  \le L_H$ ($\beta_i$ can be negative), for
    \begin{eqnarray}
        d_i > \frac{B + L_H}{6}\|s_c^{(i)}\|^{-1}>0,
        \label{dc}
    \end{eqnarray}
     we have 
    $
    m_3(s^{(i)}+s_c^{(i)}) = M(s^{(i)},s_c^{(i)}) \le M_c(s^{(i)},s_c^{(i)}). 
    $ 
\end{enumerate}
In both cases, the equality only happens at $s_c^{(i)} = 0$. 
\end{theorem}

\begin{proof}
For notational simplicity, in this proof, we use $s_c$ to represent $s_c^{(i)}$. By definition of $M_c$ and $M$, we have
\begin{eqnarray}
 M_c(s^{(i)}, s_c) -  M(s^{(i)}, s_c) =  \frac{1}{6} \bigg[ \beta_i\|s_c\|^3-T_i  [s_c]^3 \bigg]+ d_i\|s_c\|^4. 
\label{difference in Mc and M}
\end{eqnarray}
If $s_c^{(i)} = 0$, clearly $M(s^{(i)},0) = M_c(s^{(i)},0). $  Assume that $\|s_c^{(i)}\| > 0$.

\noindent
\bb{(\textbf{Case 1)}} If $\beta_i \ge L_H$, for any $d_i \ge 0$, we deduce from \eqref{difference in Mc and M} that 
\begin{eqnarray*}
    M_c(s^{(i)}, s_c) -  M(s^{(i)}, s_c )
  \ge
    \frac{1}{6} \bigg[\underbrace{L_H - \frac{T_i  [s_c]^3}{\|s_c\|^3}}_{\ge 0\text{ by Corollary \ref{corollary upper bound for Ti}}} \bigg] \|s_c\|^3 + d_i\|s_c\|^4 > 0 
\end{eqnarray*} 

\noindent
\bb{(\textbf{Case 2)}} If $\beta_i  < L_H$ ($\beta_i$ can be negative), we choose $d_i > \frac{1}{6}(B + L_H)\|s_c\|^{-1}>0$. Using 
 $ \beta_i  > -B$ from Assumption \ref{assumption1} and $T[s_c]^3 \ge -L_H \|s_c\|^3$ from Corollary \ref{corollary upper bound for Ti}, we deduce from \eqref{difference in Mc and M} that 
\begin{eqnarray*}
    M_c(s^{(i)}, s_c) -  M(s^{(i)}, s_c)> \frac{1}{6}\big(-B  - L_H \big)\|s_c\|^3 + d_i\|s_c\|^4 
    \underset{\eqref{dc}}{>} 0.
\end{eqnarray*}
\end{proof}

\begin{remark}
 If $d_i $ satisfies \eqref{dc} in Theorem \ref{thm cubic Upper Bound for s^i}, then $ \sigma_c^{(i)}$ satisfies the lower bound on $\sigma_c^{(i)} = \sigma + 4d_i$ in Algorithm \ref{algo: cqr variant 1}, as given in \eqref{addition cond on di}. This is because if $    d_i  \ge \frac{1}{6}(B + L_H)\|s_c^{(i)}\|^{-1}$, using $B > -\beta_i$ and $L_H>1>\alpha>0$, we have  
    $      d_i  \ge \frac{1}{6} \big(-\beta_i + \alpha\big)\|s_c^{(i)}\|^{-1}  $
    and, consequently,  $\sigma_c^{(i)}  = \sigma + 4d_i\ge \frac{2}{3} \big(-\beta_i + \alpha\big)\|s_c^{(i)}\|^{-1}.$
\end{remark}

{Theorem \ref{thm cubic Upper Bound for s^i} provides some intuition about how $\beta_i$ approximates the third-order information.} In particular, if $\beta_i$ is positive and the $\beta_i \|s\|^3$ term sufficiently regularizes the third-order information, then any $d_i\ge 0$ will ensure $M_c$ stays above $m_3$. On the other hand, when $\beta_i$ is negative, we need to ensure that the regularization term of $M_c$ is large enough to guarantee  {$m_3(s^{(i)}+s_c^{(i)}) \le M_c(s^{(i)},s_c^{(i)})$}. Therefore, we require $d_i$ to satisfy \eqref{dc} in Case 2.

We define the sets of successful and unsuccessful iterations for Algorithm \ref{algo: cqr variant 1} as $\mathcal{S}_i$ and $\mathcal{U}_i = \mathcal{U}_i: = [0:i]\backslash   \mathcal{S}_i$. Notice that  for $i \in \mathcal{S}_i$, $\rho_j \ge \eta$ and $s^{(i+1)} = s^{(i)}+ s_c^{(i)}$, while $s^{(i+1)} = s^{(i)}$ for  $ i \in  \mathcal{U}_i$.

\begin{lemma}
\label{lemma upper bound on di}
\textbf  {(Upper bound in $d_i$)}
Suppose that  Assumption \ref{assumption1} holds and Algorithm \ref{algo: cqr variant 1} is used. 
Then, for all $i \ge 0$
\begin{eqnarray}
 d_i \le d_{\max} :=  \gamma (B + L_H)^{3/2}\epsilon^{-1/2},
\label{uppe bound on di} 
\end{eqnarray}
where $\gamma$ is a fixed constant from Algorithm \ref{TCQR}. 
\end{lemma}

\begin{proof}
\textbf{\bb{(Case 1)}} If $\beta_i \ge L_H$, for any $d_i \in [0, d_{\max}]$, according to Theorem \ref{thm cubic Upper Bound for s^i}, we have
\begin{eqnarray*}
m_3(s^{(i)}) -m_3(s^{(i)}+s_c^{(i)}) = m_3(s^{(i)}) - M(s^{(i)}, s_c^{(i)}) \ge  m_3(s^{(i)})-   M_c(s^{(i)}, s_c) > 0. 
\end{eqnarray*}
Therefore,  $
\rho_i := \frac{m_3(s^{(i)}) -m_3(s^{(i)}+s_c^{(i)})}{m_3(s^{(i)}) - M_c(s^{(i)}, s_c^{(i)})} \ge 1\ge \eta.
$  The iteration $i$ is successful. 

\noindent
\textbf{\bb{(Case 2)}} If $\beta_i  < L_H$ ($\beta_i$ can be negative) and suppose that $$
d_i \ge  (B + L_H)^{3/2}\epsilon^{-1/2}.
$$
We deduce that $
d_i = d_i^{2/3}d_i^{1/3} \ge  (B + L_H) \epsilon^{-1/3} d_i^{1/3}. $ 
This implies the bound
\begin{eqnarray}
    d_i \ge    (B + L_H) \max\bigg\{(B + L_H)^{1/2}\epsilon^{-1/2}, {d_i}^{1/3} \epsilon^{-1/3}\bigg\} \underset{\eqref{lower bound on step size}}{\ge}  (B + L_H)\|s_c^{(i)}\|^{-1} 
    \label{dc on successful step}
\end{eqnarray}
where the last inequality follows from the upper bound on $\|s_c^{(i)}\|$ in Lemma \ref{lemma: lower bound of norms}.
Since \eqref{dc on successful step} satisfies the requirement on $d_i$ in Theorem \ref{thm cubic Upper Bound for s^i}, we deduce that $m_3(s^{(i)}+s_c^{(i)}) \le  M_c(s^{(i)}, s_c^{(i)})$. The rest of the analysis follows similarly as in Case 1. We deduce that the iteration $i$ is successful. 

In both cases, it follows from the update of regularization parameters in Algorithm \ref{algo: cqr variant 1} that the next iteration $0 \le d_{i+1} \le d_i$. As a consequence, the mechanism of regularization parameters update for Algorithm \ref{algo: cqr variant 1} ensure that \eqref{uppe bound on di} holds.
\end{proof}

\begin{remark}
It is worth noting that \eqref{addition cond on di} in Algorithm \ref{algo: cqr variant 1} is below the upper bound $d_{\max}$; {namely, there is no contradiction between \eqref{addition cond on di} and \eqref{uppe bound on di}.} To see this, assume instead that
$$
\sigma_c^{(i)} = \sigma + 4d_i \ge \frac{2}{3} \big(-\beta_i + \alpha\big)\|s_c^{(i)}\|^{-1}  \ge  \sigma + 4d_{\max}. 
$$
This implies $ \frac{1}{6} \big(-\beta_i + \alpha\big)\|s_c^{(i)}\|^{-1} \ge d_{\max} $. However, using $0< \alpha <1 \le L_H $, the lower bound for $\|s_c^{(i)}\|$ in Lemma \ref{lemma: lower bound of norms}, $d_{\max}$ in  \eqref{uppe bound on di}, $\gamma>1$ and $|\beta_i| <B$, we deduce that
\begin{eqnarray*}
    \frac{1}{6} \big(-\beta_i + \alpha\big) \|s_c^{(i)}\|^{-1} &<& \frac{1}{6} \big(B+\alpha \big)   \max\bigg\{(B  +L_H )^{1/2}\epsilon^{-1/2}, 2{d_i}^{1/3}  \epsilon^{-1/3} \bigg\}
    \\ &<&    \frac{1}{3}
  \max\bigg\{ \frac{1}{2}\underbrace{(B  +L_H )^{3/2}\epsilon^{-1/2}}_{=\gamma^{-1}d_{\max}}, \underbrace{\big(B+L_H \big)  \epsilon^{-1/3}}_{=\gamma^{-2/3}d^{2/3}_{\max}}  {d_i}^{1/3} \bigg\}  
  \\ &<&    \frac{1}{3} \max\big\{(2\gamma)^{-1}, \gamma^{-2/3} \big\}  d_{\max} < d_{\max} 
\end{eqnarray*}
which is a contraction. 
\end{remark}

The intuition behind Lemma \ref{lemma upper bound on di} is that whenever the regularization parameter exceeds the bound in \eqref{uppe bound on di}, the CQR model $M_c(s^{(i)}, s)$ overestimates $m_3$ at its minimizer $s= s^{(i)}_c$ (Theorem \ref{thm cubic Upper Bound for s^i}). Therefore, any decrease obtained in minimizing  $M_c(s^{(i)}, s)$ is a lower bound on the decrease in  $m_3$, ensuring that every iteration
is successful, and consequently removing the need to increase the regularization parameter any further.

\begin{theorem}
\label{thm M3 value decrease}  
\textbf{(A decrease in values of $m_3$ in successful iterations)}
Suppose that  Assumption \ref{assumption1} holds and Algorithm \ref{algo: cqr variant 1} is used. For $\alpha \in (0, \frac{1}{2})$ and $\eta \in (0, 1]$, assume that $i+1$ is a successful iteration {which does not meet the termination condition, i.e., $\|g_{i+1}\|> \epsilon$ and $\|g_i\| > \epsilon$}. Then,
\begin{equation}
  m_3(s^{(i)})-  m_3(s^{(i+1)})
 \ge  \frac{\alpha \eta}{24} \|s_c^{(i)}\|^3 \ge \kappa_s \epsilon^{3/2}
  \label{general model decrease}
\end{equation}
where $\kappa_s = \frac{\alpha \eta}{24}(B + L_H)^{-3/2}$. 
\end{theorem}

\begin{proof}
For notational simplicity, in this proof, we use $s_c$ to represent $s_c^{(i)}$, $\sigma_c$ to represent $\sigma_c^{(i)}$, and $d_i$ to represent $d_i^{(i)}$. Since $i \in \mathcal{S}_i$, \eqref{ratio test} gives
$$
\eta^{-1}\left[m_3(s^{(i)}) - m_3(s^{(i+1)})\right] \ge M_c(s^{(i)}, 0) - M_c(s^{(i)}, s_c) 
= - g_i^T s_c - \frac{1}{2} H_i [s_c]^2 - \frac{\beta_i}{6} \|s_c\|^3 - \frac{\sigma_c}{4}\|s_c\|^4
$$
where $s_c^{(i)}$ is obtained by Algorithm \ref{TCQR} and is the global minimum of $M_c$. Then, $s_c^{(i)}$ satisfies \eqref{cond0 nonlinear} and \eqref{cond1 nonlinear}, which implies
\begin{eqnarray} 
M_c(s^{(i)}, 0) - M_c(s^{(i)}, s_c)   \underset{\eqref{cond0 nonlinear}}{=} \frac{1}{2} H_i [s_c]^2 +  \frac{\beta_i}{3} \|s_c\|^3+ \frac{3\sigma_c}{4}\|s_c\|^4  
\underset{\eqref{cond1 nonlinear}}{\ge}   \frac{\beta_i}{12} \|s_c\|^3+ \frac{\sigma_c}{4}\|s_c\|^4.
\label{decrease on Mc}
\end{eqnarray}

\begin{itemize}
    \item \bb{A1)}: Assuming $\beta_i \in  [\alpha, B]$, since $\sigma_c >0$, we deduce from \eqref{decrease on Mc} that
$$
  m_3(s^{(i)})-  m_3(s^{(i+1)})
 \ge \frac{\beta_i\eta}{12} \|s_c\|^3+ \frac{\sigma_c\eta}{4}\|s_c\|^4>  \frac{\alpha \eta}{12} \|s_c\|^3.
$$

\item  \bb{A2)}: Assuming $ \beta_i \in  [-B, \alpha] $ and $\sigma_c \ge \frac{2}{3} \big(-\beta_i + \alpha\big)\|s_c\|^{-1} 
> 0$, we deduce from \eqref{decrease on Mc} that
\begin{eqnarray*}
  m_3(s^{(i)})-  m_3(s^{(i+1)}) \ge \frac{\beta_i \eta}{12} \|s_c\|^3+ \frac{\eta}{6} \big(-\beta_i + \alpha\big)\|s_c\|^3  
  \ge  \bigg(- \frac{\beta_i}{12} +\frac{\alpha}{6}  \bigg) \eta\|s_c\|^3 \ge \frac{\alpha \eta}{12} \|s_c\|^3. 
\end{eqnarray*}
where the last inequality follows from $-\beta_i \ge - \alpha$ for $\beta_i \in [-\alpha, \alpha] $ and $-\beta_i \ge  \alpha$ for $\beta_i \in [-B, -\alpha] $. 

\item 
\bb{A3)}: Assuming $ \beta_i \in  [-B, -4\alpha] $ and $0 < \sigma_c \le \frac{1}{6} \big(-\beta_i + \alpha\big)\|s_c\|^{-1}$: Since $s_c$ is a global minimum, $s_c$ is a local minimum. According to the second-order local optimality condition \eqref{hessm}, 
$
0 \le \nabla^2 M_c(s^{(i)}, s_c)[s_c]^2 = H_i[s_c]^2 + \beta_i \|s_c\|^3 + 3 \sigma_c \|s_c\|^4. 
$
Substituting the local optimality condition into the middle equation in \eqref{decrease on Mc}, we have
\begin{eqnarray}
M_c(s^{(i)}, 0) - M_c(s^{(i)}, s_c) \ge- \frac{\beta_i}{6}  \|s_c\|^3- \frac{3}{4} \sigma_c \|s_c\|^4.
\label{decrease on Mc 3}
\end{eqnarray}
Using $0 < \sigma_c \le \frac{1}{6} \big(-\beta_i + \alpha\big)\|s_c\|^{-1}$, we deduce from \eqref{decrease on Mc 3} that
\begin{eqnarray*}
m_3(s^{(i)}) - m_3(s^{(i+1)}) \ge - \frac{\eta\beta_i }{6} \|s_c\|^3- \frac{\eta}{8}  \big(-\beta_i + \alpha\big)\|s_c\|^3= \bigg(-\frac{\beta_i}{24} - \frac{\alpha}{8}  \bigg) \eta \|s_c\|^3 
\ge  \frac{\alpha \eta}{24} \|s_c\|^3
\end{eqnarray*}
where the last inequality uses  $ \beta_i \le -4\alpha.$ 
\end{itemize}
{
We deduce from (A1) that, in a successful step, if we have $\beta_i \in [\alpha, B]$ and $\sigma_c > 0$, then
\begin{eqnarray}
    m_3(s^{(i)})-  m_3(s^{(i+1)})  \ge   \frac{\alpha \eta}{24}\|s_c\|^3.
  \label{general model decrease temp} 
\end{eqnarray}
Also, we can deduce from (A2) and (A3) that if $\beta_i \in [-B, -4\alpha]$ and $\sigma_c \ge \frac{2}{3} \big(-\beta_i + \alpha\big)\|s_c\|^{-1}$ or $\sigma_c \le \frac{1}{6} \big(-\beta_i + \alpha\big)\|s_c\|^{-1}$, then \eqref{general model decrease temp} also holds. Lastly, using (A3), if $\beta_i \in [-4\alpha, \alpha]$ and $\sigma_c \ge \frac{2}{3} \big(-\beta_i + \alpha\big)\|s_c\|^{-1}$, then \eqref{general model decrease temp} also holds. These three cases correspond to the requirements of $\beta_i$ and $\sigma_c$ in Algorithm \ref{algo: cqr variant 1}. Using the lower bound on $\|s_c\|$ in Lemma \ref{lemma: lower bound of norms}, we deduce the desired result in \eqref{general model decrease}.}
\end{proof}

\begin{remark}
\label{remark cauchy analysis}
It is important to note that the $\sigma_c^{(i)}$ condition in Algorithm \ref{algo: cqr variant 1} are meant to ensure that the decrease in values of $m_3$ is at least $O(\epsilon^{3/2})$, specifically: $
  m_3(s^{(i)})-  m_3(s^{(i+1)})  \ge   \frac{\alpha \eta}{24}(B + L_H)^{-3/2}\epsilon^{3/2}.
$ If we only have the condition $\rho_i \ge \eta$, the CQR algorithmic framework (Algorithm \ref{TCQR}) still results in a decrease in value of $m_3$ at each successful iteration. This is evident from the sufficient condition for a global minimizer of $M_c$ in Theorem \ref{thm: cubic sufficient Nocedal}, from $\beta \ge -3\sigma_c\|s^{(i)}_c\|$, we deduce that the right-hand side of \eqref{decrease on Mc} is positive. We can also verify the decrease in value of $m_3$ with the following argument. Assume $\sigma_c > 0$; the eigenvalue of $H_i$ is bounded for all $i \ge 0$,
\begin{eqnarray*} 
  \lambda_{\max}(H_i) &=& \max_{u \in \R^n,\|u\|=1} \bigg( H + T[s^{(i)}]+ \sigma \|s^{(i)} \|^2 + 2\sigma (s^{(i)})^T s^{(i)}\bigg)[u]^2 \\ &\le&  H + \Lambda_0 \|s^{(i)} \| + 3 \sigma \|s^{(i)} \|^2 \le\|H\| + \Lambda_0  r_c + 3\sigma r_c^2 =:L_g     
\end{eqnarray*}
where $\|H\| = \lambda_{\max}[H]$ and $\Lambda_0 := \max_{\|u_1\|=\|u_2\|=\|u_3\|=1} {T [u_1][u_2][u_3]}$. 
Assume that $|\beta_i|\le B$ and set $v := -k\frac{g_i}{\|g_i\|}$ for any scalar $k>0$. 
 Since $s_c^{(i)}$ is the global minimizer, we have
\begin{eqnarray*} 
&& M_c(s^{(i)}, 0)-    M_c(s^{(i)}, s_c^{(i)})  \ge   f_i -M_c(s^{(i)}, v)
= k\|g_i\| -\frac{k^2}{2}\frac{g_i^TH_ig_i}{\|g_i\|^2}-\frac{\beta_i}{6} k^3 -  \frac{\sigma_c^i}{4}k^4 
\\  && \qquad \qquad >k\epsilon -\frac{L_g}{2}  k^2-\frac{B}{6} k^3 -  \frac{\sigma_c^i}{4}k^4= k\bigg(\frac{\epsilon}{3} -\frac{L_g}{2}  k\bigg) + k\bigg(\frac{\epsilon}{3} -\frac{B}{6} k^2 \bigg)+   k\bigg(\frac{\epsilon}{3}-\frac{\sigma_c^i}{4}k^3\bigg). 
\end{eqnarray*}
By choosing, $ 0< k < \frac{1}{2} \min \big\{ \frac{2}{3L_g}\epsilon, (\frac{B}{2}\epsilon)^{1/2}, (\frac{4}{3\sigma_c^i}\epsilon)^{1/3}\big\}$, every term in the bracket is negative, and we have  
$$
  m_3(s^{(i)})-  m_3(s^{(i+1)})  \ge \eta \bigg[ M_c(s^{(i)}, 0)-    M_c(s^{(i)}, s_c^{(i)}) \bigg]\ge \frac{\epsilon k \eta}{6} > \frac{\epsilon^2 \eta}{18L_g}. 
$$
This analysis is unaffected by $\beta_i$ being positive or negative and does not require the $\sigma_c^{(i)}$ condition in Algorithm \ref{algo: cqr variant 1}. 
\end{remark}

We may now focus on  how many successful iterations may be needed to compute an $\epsilon$-approximate first-order minimizer using Algorithm \ref{algo: cqr variant 1}.
\begin{lemma} 
\label{lemma Number of Iterations}
\textbf{(Successful and unsuccessful adaptive-regularization iterations)}
    Suppose that Algorithm \ref{algo: cqr variant 1} is used and that $d_i \le d_{\max}$. Then, 
    $$
    i \le |\mathcal{S}_i|  + \frac{\log(d_{\max})}{\log \gamma}.
    $$ 
\end{lemma}

{Lemma \ref{lemma Number of Iterations} employs a similar technique to that in \cite[Lemma 2.3.1, Lemma 2.4.1]{cartis2022evaluation}, and the proof is provided in detail in Appendix \ref{appendix: proof for lemma 3.6}.} Note that Lemma \ref{lemma Number of Iterations} is entirely independent of the form of the model $M_c$ and depends on the mechanism defined by the update of regularization parameters in Algorithm \ref{algo: cqr variant 1}. We do not need  Assumption \ref{assumption1} for Lemma \ref{lemma Number of Iterations}.

\begin{lemma}
\label{lemma Bound on the number of successful iterations}
\textbf{(Bound on the number of successful iterations)} Suppose that  Assumption \ref{assumption1}  holds.  Let $m_{\text{low}}$ be a lower bound on $m_3(s)$ for $s \in \R^n$. Then, there exists a positive constant $\kappa_s$  as defined in Theorem \ref{thm M3 value decrease} such that Algorithm \ref{algo: cqr variant 1} requires at most
\begin{eqnarray}
  |\mathcal{S}_i|:= \kappa_s \frac{m_3(s^{(0)})-m_{\text{low}}}{\epsilon^{3/2}}
  \label{bound for successful iterations}
\end{eqnarray}
successful iterations before an iterate $s_\epsilon$ is computed for which $\|g_i\| \le \epsilon$.
\end{lemma}
\begin{remark}
    Given $\sigma>0$ as defined in \eqref{m3},  $m_3(s)$ is bounded below for $s \in \R^n$.
\end{remark}
\begin{proof}
Suppose that $\|g_j\| > \epsilon$ {and $\|g_{j+1}\| > \epsilon$} for all $j \in [1: i]$. Then, 
\begin{eqnarray*}
m_3(s^{(0)})-m_{\text{low}} \ge m_3(s^{(0)})-m_3(s^{(i+1)}) = \sum_{j \in \mathcal{S}_i} \bigg[m_3(s^{(i)})-m_3(s^{(i+1)}) \bigg] >  |\mathcal{S}_i|\kappa_s \epsilon^{3/2}
\end{eqnarray*}
where $\kappa_s$ given in Theorem \ref{thm M3 value decrease}. This implies the bound \eqref{bound for successful iterations}. 
\end{proof}

Given Lemma \ref{lemma Bound on the number of successful iterations}, we are now ready to state the worst-case evaluation bound for Algorithm \ref{algo: cqr variant 1}. 

\begin{theorem} \textbf{(Complexity bound for Algorithm \ref{algo: cqr variant 1})} Suppose that  Assumption \ref{assumption1}  holds.  Then, there exists a positive constant $\kappa_s$ (see Theorem \ref{thm M3 value decrease}),  $\kappa_c$ such that Algorithm \ref{algo: cqr variant 1} requires at most
$$
 \kappa_s \frac{m_3(s^{(0)})-m_{\text{low}}}{\epsilon^{3/2}}  + \kappa_c + \log{\epsilon^{-1/2}} = O(\epsilon^{-3/2}). 
$$
function evaluations (i.e., $f_i$) and at most
$$
 \kappa_s \frac{m_3(s^{(0)})-m_{\text{low}}}{\epsilon^{3/2}}  + 1
= O(\epsilon^{-3/2}). 
$$
derivatives evaluations (i.e., $g_i, H_i, T_i$) to compute an iterate $s_\epsilon$ such that $\|g_i\| \le \epsilon$ {or $\|g_{i+1}\| \le \epsilon$} .
\label{thm: cqr complexity}
\end{theorem}

\begin{proof}
    The number of successful iterations needed to find an $\epsilon$-approximate first-order minimizer is bounded above by \eqref{bound for successful iterations}. The total number of iterations (including unsuccessful ones) can also be bounded above using Lemma \ref{lemma Number of Iterations}, which yields that this number is at most
$$
 \kappa_s \frac{m_3(s^{(0)})-m_{\text{low}}}{\epsilon^{3/2}}  + \frac{\log(d_{\max})}{\log \gamma} \underset{\eqref{uppe bound on di}} {=}  \kappa_s \frac{m_3(s^{(0)})-m_{\text{low}}}{\epsilon^{3/2}}  + 1 + \frac{3 \log (B + L_H)}{2\log\gamma} + \log{\epsilon^{-1/2}}. 
$$
where we use $d_{\max} =  \gamma (B + L_H)^{3/2}\epsilon^{-1/2}$. 
Since Algorithm \ref{algo: cqr variant 1} uses at most one evaluation of $m_3$ per iteration and at most one evaluation of its derivatives per successful iteration (plus one evaluation of $m_3$ and its derivatives at the final iteration), by setting $\kappa_c :=  1  + \frac{3 \log(B + L_H)}{2\log\gamma}$, we deduce the desired conclusion. 
\end{proof}

\subsection{CQR Algorithm in Some Special Cases}
\label{sec: CQR Algorithm in Special Cases}

\begin{figure}[ht]
    \centering
    \includegraphics[width = 8cm]{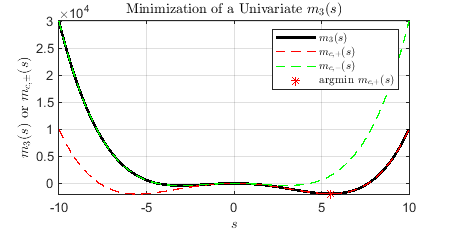}
    \caption{\small Minimization of univariate $m_3(s) = 3s -50s^2 -10 s^3 + 3s^4$. }
    \label{univariate m}
\end{figure}

\textbf{Univariate $m_3$: convergence of the CQR algorithm in one step.} For univariate cases ($n=1$), we define $W$ as $W := 1$, $\sigma_c^{(0)}:=\sigma$, and $\beta^{(0)} := \pm |T| \in \mathbb{R}$. We introduce two univariate CQR polynomials, denoted as $m_{c, +}$ and $m_{c, -}$, as follows:
\begin{align*}
    m_{c, +}(s) := f_0 + gs + \frac{1}{2} Hs^2 + \frac{|T|}{6}|s|^3 + \frac{\sigma}{4}s^4, \qquad
    m_{c, -}(s) := f_0 + gs + \frac{1}{2} Hs^2 - \frac{|T|}{6}|s|^3 + \frac{\sigma}{4}s^4.
\end{align*}
As depicted in Figure \ref{univariate m}, $m_{c, -}$ or $m_{c, +}$ overlap with $m_3$ in one of the half-spaces $s>0$ or $s<0$. Therefore, to determine the global minimum of the univariate function $m_3$, our objective is to minimize the CQR polynomials $m_{c, -}$ and $m_{c, +}$. This involves finding $s_{c, +}=\argmin_{s \in \mathbb{R}} m_{c, +}(s) $ and $ s_{c, -}=\argmin_{s \in \mathbb{R}} m_{c, -}(s)$
using Algorithm \ref{Newton algo}. The minimizer of $m_3$ is denoted as $s^* = \argmin_{s \in \R} m_3(s) = \argmin \big\{m_3( s_{c, +}), m_3(s_{c, -}) \big\}$. We can determine the optimal value of the univariate function $m_3$ within a single iteration of the CQR algorithm (i.e. minimizing CQR polynomials $m_{c, -}$ and $m_{c, +}$ once each). An illustrative example of this minimization process is presented in Figure \ref{univariate m}, where the global minimizer for $m_3(s)$ aligns with that of $m_{c, +}(s)$.

In the multivariate case ($n \ge 2$), the tensor directions are not restricted to $s \ge 0$ or $s \le 0$ as they are in the univariate case. Therefore, it is generally not possible to find a choice of $m_{c, \pm}(s)$ that is equal to $m_3(s)$. However, similar to the univariate case, we can formulate the half plane as $\mathcal{H}_+ = \{s_c^{(i)} \in \R^n,  T[s_c^{(i)}]^3 \ge 0\}$ and $\mathcal{H}_- = \{s_c^{(i)} \in \R^n,  T[s_c^{(i)}]^3 < 0\}$. Our aim is to choose {$\beta_i \approx  \frac{T[s_c^{(i)}]^3}{\|s_c^{(i)}\|}$. More details can be found in Section \ref{sec cqr numerical Algorithm}.}  

\textbf{Diagonal Hessian and Tensor:} Suppose that $H = \nabla^2 f(x_k)$ and $T = \nabla^3 f(x_k)$ are diagonal matrices and tensors with diagonal entries $D_\iota \in \R$ and $\tilde{T}_\iota \in \R$ for $\iota = 1, \dotsc, n$, respectively. By incorporating a slightly different fourth-order regularization {norm, namely $ \sigma \| s\|_4^4=  \sigma\sum_{i=1}^n s_\iota^4$ with $\sigma>0$}, we can obtain a separable cubic-quartic polynomial denoted as
{$$
\tilde{m}_3(s) : =   f_0 +  g^T s+ \frac{1}{2}H[s]^2 + \frac{1}{6}T[s]^3 + \sigma \| s\|_4^4 = f_0 +  \sum_{\iota=1}^n \bigg[ \underbrace{g_\iota s_\iota+\frac{1}{2} D_\iota s_\iota^2  +  \frac{1}{6} \tilde{T}_\iota s_\iota^3 + \frac{\sigma}{4} s_\iota^4}_{\text{univariate CQR Polynomial}} \bigg].
$$}
Here, $s_\iota \in \R$ and $g_\iota \in \R$ denote the $\iota$th entry of the vectors $s \in \R^n$ and $g \in \R^n$, respectively.
Consequently, we expressed $\tilde{m}_3(s)$ as the sum of $n$ univariate CQR polynomials, and minimizing $\tilde{m}_3(s)$ can be reduced to the minimization of $n$ separate univariate CQR polynomials. Each of these polynomials can be efficiently solved using a single step of the univariate CQR algorithm. 

\textbf{Small Tensor Term: }If the tensor term is small, $\|T_i\| \le \epsilon^{1/3}$, and we also choose $|\beta_i| \le \epsilon^{1/3}$, then \eqref{diff grad} becomes
$
\epsilon< \frac{\epsilon^{1/3}}{2}\|s_c\|^2 +    \frac{\epsilon^{1/3}}{2} \big\|s_c\big\|^2+ 4d_i\|s_c\|^3.
$
This inequality implies $\|s_c^{(i)}\| > O(\epsilon^{1/3})$. Following the analysis in Theorem \ref{thm M3 value decrease}, this results in a complexity bound of order $O(\epsilon)$. 

\section{Numerical Implementation and Preliminary Results}
\label{sec cqr numerical Algorithm}
In this section, we present the implementation of CQR method and provide preliminary numerical results of the algorithm. The CQR Variant for numerical implementation is described in Algorithm \ref{PCQR}. 

\begin{algorithm}
Initialize and proceed through Steps 1 to 4 of the CQR algorithmic framework (Algorithm \ref{TCQR}). Specifically, in each iteration $(i)$, compute ${s_c^{(i)}} := \argmin_{s\in \R^n} M_c(s^{(i)}, s)$, where $\beta_i$ is initialised and updated using one of the following updating rules: 
\begin{enumerate}
    \item $
\beta_{i+1} := L_H(s^{(i)}) = \Lambda_0 + 6 \sigma \|s^{(i)}\|$ for $i \ge 0$; $\qquad \beta_0 := \Lambda_0 = \sqrt{\sum_{\iota}\max_{1 \le  j, k \le n} T[ \iota, j, k]^2}$. 
    \item $
\beta_{i+1} := \frac{T_i [s_c^{(i)}]^3}{\|s_c^{(i)}\|^{3}}$ for $i \ge 0$; $\qquad \qquad \qquad \qquad  \qquad\beta_0 := - \max_{1 \le j \le n} |T[j, j, j]|$. 
    \item  $
\beta_{i+1}  := \frac{1}{n}\sum_{j=i}^n T_i[j, j, j]$ for $i \ge 0$; $\qquad\qquad \qquad \beta_0 := - \max_{1 \le j \le n} |T[j, j, j]|$. 
\end{enumerate}
\textbf{In choices 2 or 3 above}: initialise an upper bound $B>0$ at the start of the algorithm; at iteration $(i)$,
if $|\beta_{i+1}| >B$, then set\footnotesize{$^b$} \normalsize{} $\beta_{i+1} :=\sign(\beta_{i+1})B$.

\footnotesize{$^b$See Remark \ref{remark bound B}.}
\caption{\small A CQR Variant for Numerical Implementation}
\label{PCQR}
\end{algorithm}

When comparing the CQR Variant for numerical implementation (Algorithm \ref{PCQR}) to the CQR algorithmic framework (Algorithm \ref{TCQR}), it is worth noting that Algorithm \ref{TCQR} does not specify a particular choice for $\beta_i$, whereas Algorithm \ref{PCQR} offers three options for updating $\beta_i$ that we can choose from. The three choices are summarised as follows. 

\begin{figure}[!ht]
    \centering
        \caption{\small Performance profile of Algorithm \ref{PCQR} using the first choice of $\beta_i$. The detailed setup of the test problem is in Section \ref{sec numerical setup} and parameters in \eqref{numerical setu} are set at \texttt{a} $=80$, \texttt{b} $=80$, \texttt{c} $=80$, $\sigma =5$, and $n=100$. }
    \label{fig beta choice 1}
    \includegraphics[width =15.5cm]{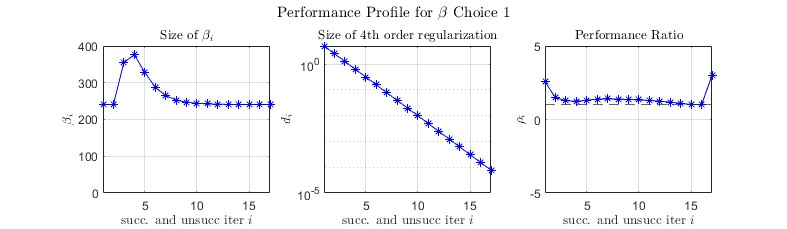}
\end{figure}

The first choice of $\beta_i$ in Algorithm \ref{PCQR} stems from our theoretical analysis in Theorem \ref{thm cubic Upper Bound for s^i} and Corollary \ref{corollary upper bound for Ti}. By setting $\beta_i = L_H(s^{(i)}) = \Lambda_0 + 6 \sigma \|s^{(i)}\|$ where {$ \Lambda_0$ is calculated as in \cite[Thm 2.6]{cartis2015branching}},  we ensure that $M(s^{(i)},s_c^{(i)}) \le M_c(s^{(i)},s_c^{(i)})$ for all $i$. This guarantees that $\rho_i>1$ in every iteration of Algorithm \ref{PCQR}, namely every step is a successful step.  
A numerical example of this choice of $\beta_i$ is given in Figure \ref{fig beta choice 1} where we plot the value of $\beta_i$, the size of $d_i$, and the performance ratio for both successful and unsuccessful iterations. Our experimental findings affirm our theoretical results that each iteration indeed satisfies $\rho_i \ge 1$ (right plot of Figure \ref{fig beta choice 1}). Consequently, there are no unsuccessful iterations. Furthermore, as each iteration achieves $\rho_i \ge 1$, which is considered as `a very successful step', we observe a reduction in $d_i$ by a factor of $\gamma_2$ per iteration (the second plot of Figure \ref{fig beta choice 1}). In the last few iterations, $d_i$ has magnitude $10^{-3}$, nearly reaching zero, while the iterations remain successful.  

The second and third choices of $\beta_i$ in Algorithm \ref{PCQR} are the practical ones for efficient algorithm implementation. In subsequent sections, we conduct numerical tests using one of these two choices. The second choice of $\beta_i$ in Algorithm \ref{PCQR} is derived from the observation in \eqref{difference in Mc and M}. For each iteration, Theorem \ref{thm cubic Upper Bound for s^i} holds if the condition $\beta_{i+1} \ge T_i [s_c^{(i+1)}]^3\|s_c^{(i+1)}\|^{-3}$ is satisfied. However, prior to minimizing $M_c(s^{(i)}, s)$, $s_c^{(i)}$ remains unknown. Therefore, in the algorithmic implementation, we update the $
\beta_{i+1}:= T_i [s_c^{(i)}]^3\|s_c^{(i)}\|^{-3}$ which is the tensor applied to the search direction from the last successful iteration. We use this update as an estimate of $T_i [s_c^{(i+1)}]^3\|s_c^{(i+1)}\|^{-3}$.  For the third choice of $\beta_i$ in Algorithm \ref{PCQR}, we opt for $\beta_i$ to be the scaled trace of $T_i$: $\beta_{i+1} := \frac{1}{n}\sum_{j=i}^n T_i[j, j, j]$. 
Preliminary examples for these two choices of $\beta_i$ are provided in Appendix \ref{appendix beta choice}.

\begin{figure}[!ht]
    \centering
   \caption{\small  Performance comparison between Algorithm \ref{algo: cqr variant 1} and \ref{PCQR} with the second choice of $\beta_i$ and $\alpha =0.1$: Algorithm \ref{algo: cqr variant 1} uses $14$ total iterations while Algorithm \ref{PCQR} uses $13$ total iterations. The differences in iterations between the two algorithms are with a red asterisk \rr{$\ast$}. The same test example is used as in Figure \ref{fig beta choice 1}.}
    \label{fig algo 2 and 3}
    \includegraphics[width =15.5cm]{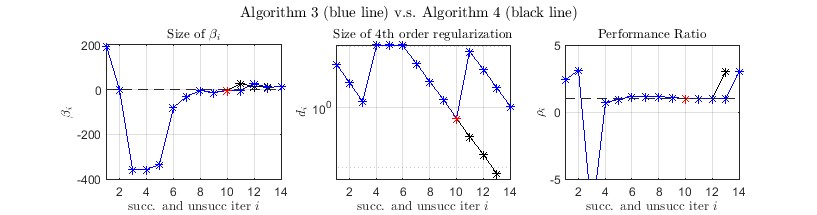}
\end{figure}

When contrasting the CQR Variant for numerical implementation (Algorithm \ref{PCQR}) with the CQR variant featuring optimal complexity (Algorithm \ref{algo: cqr variant 1}), a notable distinction arises. Unlike the Algorithm \ref{algo: cqr variant 1}, Algorithm \ref{PCQR} does not have the additional requirement on $\sigma_c^{(i)}$, thereby allowing for a more adaptive choice of $\sigma_c^{(i)}$. In other words, Algorithm \ref{algo: cqr variant 1} adopts a more conservative approach when accepting steps compared to Algorithm \ref{PCQR}. Yet, according to the analysis presented in Remark \ref{remark cauchy analysis}, Algorithm \ref{PCQR} ensures a reduction in the value of $m_3$ during each successful iteration and converges to a point satisfying $\|g_i\| < \epsilon$. A numerical example comparing Algorithm \ref{algo: cqr variant 1} and Algorithm \ref{PCQR} is given in Figure \ref{fig algo 2 and 3}.

\begin{remark}
{In Assumption \ref{assumption1}, we assume that $\beta_i$ is bounded for all iterations. In Algorithm \ref{PCQR}, the parameter $B$ prevents $|\beta_i|$ from becoming excessively large, thereby averting potential numerical instability.}
\label{remark bound B}
\end{remark}



{Choices 1--3 for updating $\beta_i$ in Algorithm \ref{PCQR} are just a few options within the general CQR framework. Other methods of selecting $\beta_i$ are also possible. For example, one could use  $\beta_i$ as the smallest or largest eigenvalue of the tensor, the smallest or largest entries of $T$, or update $\beta_i$ in an unsuccessful step. Further investigation into these $\beta_i$ choices will be pursued in future research.}

\subsection{Numerical Results: Comparison Across Methods}
\label{sec: numerical testing}

In this section, we carry out numerical experiments for CQR algorithm (Algorithm \ref{PCQR}) with practical $\beta$ choices (Choices 2 or 3). To assess the performance of CQR solvers, we devise multiple sets of AR$3$ subproblems across different dimensions, $n$, ranging from $5$ to $100$. These test sets comprise AR$3$ subproblems with varying gradients, Hessians, tensors, and regularization terms. Specifically, we consider singular, ill-conditioned, indefinite, and diagonal Hessians, as well as large and small tensor terms, diagonal tensor terms, badly scaled tensor terms, and dense tensor terms. We also provide examples arising from the minimization of multi-dimensional objective functions that illustrate further aspects of CQR method's performance.

\subsubsection{Numerical Set-up}
\label{sec numerical setup}
The test sets are constructed as follows: we generate quartic regularized polynomials to assess the algorithm's performance. Specifically, we define $m_3(s) = f_0 + g^T[s] + \frac{1}{2} H [s]^2 + \frac{1}{6} T [s]^3+ \frac{1}{4} \sigma \|s\|_2^4$, where $f_0=0$, and the coefficients $g$, $H$, and $T$ are generated as follows:
\begin{eqnarray}
\quad g = \texttt{a*randn(n, 1)}, \quad  H = \texttt{b*symm(randn(n, n))},  \quad  T = \texttt{c*symm(randn(n, n, n))}
\label{numerical setu}
\end{eqnarray}
Here, \texttt{n} represents the dimension of the problem ($s\in \R^n$), and \texttt{symm(randn())} denotes a symmetric matrix or tensor with entries following a standard normal distribution with mean zero and variance one. The parameters \texttt{a}, \texttt{b}, \texttt{c}, and $\sigma$ are selected differently to test the algorithm's performance under various scenarios. We set the stopping criterion to be $\|g_i\| < \epsilon$, where $\epsilon= 10^{-5}$ unless otherwise specified.
The parameters in Algorithm \ref{TCQR} are set as $\eta_1 = 0.9$, $\eta = 0.1$, $\gamma_2 = 0.5$, $\gamma = 2.0$ and $B = \max_{1 \le \iota, j, k \le n} \big|T[ \iota, j, k]\big|$. The iteration counts include both successful and unsuccessful iterations. The counts of function and derivative evaluations match the number of successful iterations. {It is worth noting that the CQR algorithm (Algorithm \ref{PCQR}) is tensor-free and does not require information about the tensor entries or its structure. Instead, the algorithm only requires tensor-vector and tensor-matrix products. }

\subsubsection{Comparison Across Methods}
In this section, we give a comparison of several minimization algorithms for the AR$3$ subproblem. These algorithms encompass the ARC method, Nesterov's method (specifically for convex $m_3$), \texttt{QQR-v1} method, and \texttt{QQR-v2} method. The ARC method, as proposed in Cartis et al. \cite{cartis2011adaptive}, forms the basis of this comparison. Nesterov's method, introduced in \cite{Nesterov2022quartic}, is tailored for convex subproblems only. The QQR method, a recent development introduced in \cite{cartis2023second}, can be considered as a nonconvex generalization of Nesterov's method; it approximates the third-order tensor term through a linear combination of quadratic and quartic terms. In our numerical testing, we assess both variants of QQR (\texttt{QQR-v1} and \texttt{QQR-v2}). \texttt{QQR-v1} follows a framework similar to Nesterov's method, employing a single adaptive parameter to govern the quadratic and quartic terms. In contrast, \texttt{QQR-v2} employs two adaptive parameters to manage the convexity scenario and regularization magnitude, respectively.

The numerical results (Tables \ref{table convex H}--\ref{table large tensor} in Appendix \ref{appendix: Numerical Results on Various Test Sets}) consistently demonstrate that our proposed methods, the CQR method (with $\beta$ choices 2 or 3 in Algorithm \ref{PCQR}), perform competitively with the ARC method and the Nesterov's method. In particular, CQR methods typically require fewer evaluations of functions/derivatives or iterations than the ARC method, Nesterov's method, and \texttt{QQR-v1}. Moreover, CQR methods exhibit robustness when handling ill-conditioned Hessian terms and perform effectively for minimizing $m_3$ with singular or diagonal Hessians (Table \ref{table ill con H}). In the case of convex and locally strictly convex $m_3$ functions  (Table \ref{table convex H}), the CQR algorithm (Algorithm \ref{PCQR}) terminates within only a few iterations and function evaluations. Notably, in these scenarios, the CQR method requires fewer iterations and evaluations to converge than Nesterov's method (which is specifically for convex $m_3$). The primary reason is that our implementation of the CQR method employs $\beta_i$ to estimate the tensor direction, whereas Nesterov's method (and \texttt{QQR-v1}) solely utilizes a single-parameter linear combination of the local upper and lower bounds.
Furthermore, the CQR method effectively locates a point satisfying $\|g_i\| < \epsilon$ for quartic polynomials charactersed by diagonal tensor terms, and large, negative, directional tensors (Tables \ref{table small tensor}--\ref{table large tensor}). Lastly, our proposed algorithm successfully identifies a point satisfying $\|g_i\| < \epsilon$ across a wide range of positive $\sigma$ values (Table \ref{change sigma}), allowing for dynamic $\sigma$ adjustments while employing CQR method to solve the AR$3$ problem. 

\subsubsection{CQR Performance in Challenging Scenarios}
\label{sec special cqr cases}
As observed in Table \ref{table convex H} to Table \ref{table large tensor}, for well-scaled subproblems, CQR methods tend to perform comparably with the \texttt{QQR-v2} method in terms of iterations and function evaluations.
However, we have noticed that in certain practical problems, such as Problem 4 of Moré, Garbow, and Hillstrom’s test set \cite{more1981testing}, known as the Brown badly scaled problem, our investigation revealed a significant performance advantage of CQR when compared to ARC and QQR. 
{We solved this problem using the AR$3$ method as described in \cite{cartis2020concise} by minimizing \eqref{ar3 model} with various subproblem solvers and compared their performances. The performance comparison is provided in Table \ref{table more 4}.} We observe that the AR$3$ method with the CQR subproblem solver requires only approximately $2$ iterations (or evaluations) per subproblem to converge, while the AR$3$ method in \eqref{ar3 model} with the ARC or QQR subproblem solver requires at least three times the iterations (or evaluations). 
After examining the subproblem linked to this challenging scenario, we noticed that the gradient, Hessian, and tensor exhibited significantly larger magnitudes in a specific entry (or in a few entries), while all other entries remained relatively small, typically around a magnitude of 1. Furthermore, the search directions were usually predominantly influenced by these tensor directions. 

\begin{table}[!htbp]
\centering
\caption{\small The objective function to minimize is Problem 4 of Moré, Garbow, and Hillstrom’s test set, known as the Brown badly scaled problem. The first-order optimality tolerances for minimizing the objective function and the subproblem are $10^{-4}$ and $10^{-5}$ respectively, and the initialization is $x_0 = \boldsymbol{0}$ with $\sigma = 1$. The AR$3$ model is outlined in \eqref{ar3 model}, more details can be found in \cite{cartis2020concise}.}
    \label{table more 4}
{\small
\begin{tblr}{
  colspec = {c|rr},
  cell{1}{2,4} = {c=2}{c},
}
\hline[2pt]
& \quad No. of iterations per subproblem &\SetCell[c=2]{}Iterations.&  \quad & \\
\hline[1pt]
Method  &  Total Iter.  &  Successful Iter.&\\
\hline[0.5pt]
\textbf{\color{blue}AR$3$ with ARC} & 16.4& 6.3 &\\
\textbf{\color{blue}AR$3$ with QQR-v2}    & 10.3& 10.2 &\\
\textbf{\color{red}AR$3$ with CQR-Choice 2}   & 1.9& 1.9 & \\
\textbf{\color{red}AR$3$ with CQR-Choice 3}   & 1.9 &1.9& \\
\hline[2pt]
\end{tblr}
}
\end{table}

Given these challenging conditions, we conducted a series of tests to assess the performance of the CQR solver in handling such badly scaled problems. 
We first design a test scenario where all other entries in $g$, $H$, and $T$ were set to 0, with only the first diagonal entry varying in magnitude from $0.1$ to $10^{5}$.  In these examples, as illustrated in the right two plots of Figure \ref{fig special example 0} in Appendix \ref{appendix: Special Examples}, CQR converges in just 1 or 2 iterations, whereas ARC and \texttt{QQR-v2} required 6 to 15 iterations. It is noteworthy that in this problem, the CQR method is also relatively insensitive to the problem size. Specifically, the CQR method converges within a few iterations throughout all dimensions. 
We attribute this behaviour to the problem's strong dependence on the direction charactersed by the largest tensor entry. The successful identification of a point that satisfies $\|g_i\| < \epsilon$ was closely linked to how effectively the algorithms progressed along the dominant tensor direction.  The CQR method excelled in efficiency in these cases.

\begin{figure}[!ht]
    \centering
\caption{\small Performance profile of a badly scaled problem: the first entry of $g$, $H$, and $T$ scales from $-0.1$ to $-10^5$, while the rest of the entries follow a normal distribution $\mathcal{N}(0, 0.1)$, and $\sigma = 50$. The left two plots have dimension $n=2$, and the right two plots have dimension $n=100$. The tensor scale is defined as the maximum absolute value of the tensor entries divided by the mean absolute value of tensor entries. $\beta_i$ is chosen as per choice 2  in Algorithm \ref{PCQR}), with similar performance observed for $\beta_i$ as per choice 3 (Figure \ref{choice 3} in Appendix \ref{appendix: Special Examples}).}
    \label{bad scale 0.1}
    \includegraphics[width =7.5cm]{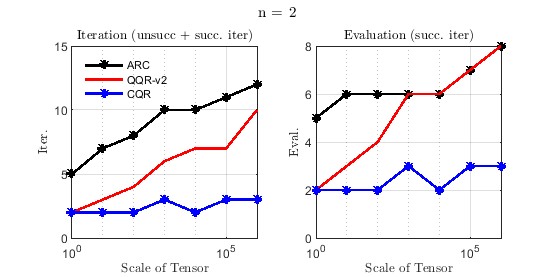}
        \includegraphics[width =7.5cm]{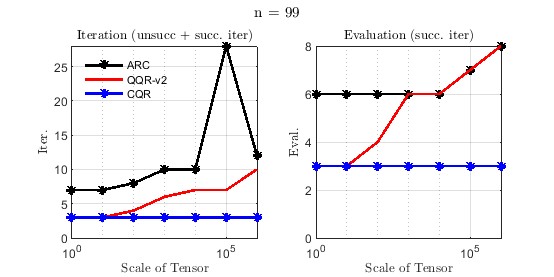}
\end{figure}

We could also add random perturbations to the zero entries of the  $g$, $H$, and $T$. We choose the gradient, Hessian, and tensors with entries of relatively small magnitudes, randomly following a normal distribution with mean zero and variance $0.1$, while the first entry of $g$, $H$, and $T$ was systematically varied from small to large magnitudes in the same direction. We tested these scenarios on subproblems ranging in size from $n=2$ to $n=100$. 
{As illustrated in Figure \ref{bad scale 0.1}, when the tensor term is well-scaled, \text{QQR-v2} and the CQR method perform similarly.} However, as the magnitude of the first entry increases, signifying a shift in the descent direction toward dominance by the tensor direction $T$, the CQR method starts to outperform ARC and \texttt{QQR-v2} in terms of both iteration  and evaluation counts. 
In particular, for a tensor scaling of $10^6$, the CQR algorithm (Algorithm \ref{PCQR}) terminates in just 4 iterations, while ARC and \texttt{QQR-v2} require 10 to 25 iterations.  
This improvement could be primarily attributed to the $\beta$ parameter, which effectively captures the dominant tensor direction that exerts a substantial influence on the subproblem.

In our exploration of additional examples, we also considered scenarios where the first diagonal entry was large, and we introduced variations in some off-diagonal entries of the tensor $T$ (Figure \ref{fig special example diag} Appendix \ref{appendix: Special Examples}). As expected, similar trends are observed, as the descent direction of $m_3$ remains predominantly influenced by the diagonal tensor direction $T$. When we expanded our experiments to include more entries with large magnitudes, the performance of CQR began to {match} that of \texttt{QQR-v2} and ARC. This behaviour was attributed to the introduction of multiple large diagonal entries in $T$, resulting in a higher rank tensor. The CQR method, with its single $\beta$ parameter, faces limitations in capturing all the tensor directions introduced.


\section{Conclusions and Future Work}
\label{sec Conclusion}

In conclusion, this paper introduces the CQR method, designed for the efficient minimization of nonconvex quartically-regularized cubic polynomials, such as the AR$3$ subproblem. The CQR method iteratively minimizes a local quadratic model incorporating cubic and quartic terms, known as the CQR polynomial. The cubic term has a coefficient $\beta_i$ that provides crucial information about the local tensor $T_i[s]^3$, while the quartic term has a coefficient $\sigma_c$ that controls the regularization and algorithmic progress. We have derived necessary and sufficient global optimality conditions for the global minimizer of the CQR polynomial. These conditions allow us to transform the problem of finding the global minimizer of the CQR polynomial into a more manageable nonlinear eigenvalue problem. Future research may explore the convergence of the root-finding algorithm for this system and scalable techniques for solving it.

Our theoretical analysis affirmed the convergence of the CQR method, showing that it is guaranteed to find a first-order critical point for the AR$3$ subproblem within at most  $\mathcal{O}(\epsilon^{-3/2})$ function and derivative evaluations; here, $\epsilon$ represents the prescribed first-order optimality tolerance. The complexity bound for the CQR method is as least as good as that of ARC, and in specific cases, it exhibits improved convergence behaviour. 
Preliminary numerical results also suggest competitive performance by our tensor-free CQR method when compared to state-of-the-art approaches such as the ARC method, often requiring fewer iterations or evaluations. Notably, our findings particularly highlight the CQR method's performance on badly scaled problems, where it outperformed other approaches like ARC and QQR in terms of both iteration and evaluation counts. 

The proposed CQR methods offer a way to approximate the third-order tensor term by a single cubic term. In the future, we plan to incorporate additional tensor information into cubic terms in CQR-like models, to capture third-order directions more accurately. By proposing efficient AR$3$ minimization algorithms, this paper brings high-order tensor methods closer to practical applications.


\appendix


\section{Derivative for $\Psi(\lambda)$}
\label{derivative appendix}
Let $H_W(\lambda) = H_i+ \lambda W$ and $W^{-1} H_W(\lambda)$ has the Cholesky factorization 
\begin{eqnarray}
W^{-1} H_W(\lambda) = L(\lambda) L^T(\lambda).
\label{chol appendix}
\end{eqnarray}
Differentiate $(H_i+ \lambda W) s(\lambda) = -g_i$ against $\lambda$ gives $    (H_i+ \lambda W)  \nabla_\lambda s(\lambda) + W s(\lambda) =0. $ Therefore, 
\begin{eqnarray}
\nabla_\lambda  s(\lambda) = -H_W^{-1}(\lambda) Ws(\lambda).
\label{ds appendix}
\end{eqnarray}
On the other hand, $\Psi(\lambda) = \|s(\lambda)\|^2_W$ has derivative $\Psi'(\lambda)  =  2 \big\langle Ws(\lambda), \nabla_\lambda s(\lambda)\big\rangle. $ Substituting \eqref{chol appendix} and \eqref{ds appendix} gives, 
\begin{eqnarray*}
\Psi'(\lambda)  &=& 
  - 2  \big\langle Ws(\lambda),  L(\lambda)^{-T}L^{-1}(\lambda) s(\lambda) \big\rangle 
\\  &=& -2  \big\langle L(\lambda)^{-1} W s(\lambda),  L(\lambda)^{-1} s(\lambda) \big\rangle = -2\|L(\lambda)^{-1}  s\|_W = -2 \|\omega(\lambda)\|_W^2. 
\end{eqnarray*}
where $L(\lambda) \omega (\lambda) = s(\lambda)$.

\section{Uniform Upper Bound for Iterates}
\label{sec: bound on iterations}

In this section, we prove that the iterates of Algorithm \ref{TCQR}, including both successful and unsuccessful iterations, are uniformly bounded above by a constant $r_c$ independent of the iteration count. The constant $r_c$ is determined by the coefficients in $m_3$, specifically $f_0$, $g$, $H$, and $T$.  We first require a technical lemma.

\begin{lemma}
Let $M_c$ be defined as in \eqref{t mc} in Algorithm \ref{TCQR} with $d_i\ge 0$. Denote $\hat{s} := s^{(i)} + s$. If $\|\hat{s}\| = \|s^{(i)}+ s\| \ge \|s^{(i)}\|$, then
    $$
M_c(s^{(i)}, s)  \ge \underbrace{f_0 + g^T[s^{(i)}+ s] +\frac{1}{2}H [s^{(i)}+ s]^2 + \frac{1}{6} T[s^{(i)}+ s]^3 }_{=T_3(\hat{s})} -  \frac{1}{6}T[ s]^3  + \frac{\beta_i}{6}\|s\|^3 + \frac{ \sigma}{20} \big \| s^{(i)}+ s\big\|^4.
    $$
\label{technical lemma}
\end{lemma}

\begin{proof}
 Using the condition $\| s^{(i)}+ s\|  \ge \|s^{(i)}\|$, we obtain that
\begin{eqnarray}
\| s^{(i)}\|^2 + \| s\|^2 + 2  s^Ts^{(i)} \ge  \|s^{(i)}\|^2  \qquad \Rightarrow \qquad  s^Ts^{(i)}  + \frac{1}{2}\|s\|^2 \ge 0. 
\label{norm eq}
\end{eqnarray}
Also, by Taylor expansion at $s=s^{(0)} =0$, 
\begin{eqnarray*}
    f_i &=&m_3(s^{(i)}) = f_0 + g^Ts^{(i)} +\frac{1}{2}H [s^{(i)}]^2 + \frac{1}{6} T[s^{(i)} ]^3 + \frac{\sigma}{4}\|s^{(i)} \|^4, 
    \label{fi}
    \\ g_i^T s  &=& \nabla m_3(s^{(i)}) =g^T s + H[s^{(i)}] [s]  + \frac{1}{2} T[s^{(i)}]^2[s] + \sigma \|s^{(i)} \|^2 s^Ts^{(i)} , 
        \label{gi}
    \\ H_i [s]^2 &=& \nabla^2 m_3(s^{(i)}) =  H[s]^2 + T[s^{(i)}][s]^2+ \sigma \|s^{(i)} \|^2 \|s \|^2 + 2 \sigma [s^T s^{(i)}]^2.
   \label{Hi}
\end{eqnarray*}
Substituting the above expressions into $M_c(s^{(i)}, s)  =  f_i + g_i^Ts+\frac{1}{2}H[s_i]^2+\frac{\beta}{6}\|s\|^3+\frac{\sigma_c}{4}\|s\|^4$ and rearranging, we have
\small
\begin{eqnarray*}
  M_c(s^{(i)}, s)  = f_0 + g^T(s^{(i)}+s) + \frac{1}{2}\bigg(H [s^{(i)}]^2 + 2 H [s^{(i)}] [s] + H [s]^2\bigg) + 
\\
  + \frac{1}{6}\bigg(  T[s^{(i)}]^3+ 3 T[s^{(i)}]^2[s] +  3T[s^{(i)}][s]^2 + \beta_i \|s\|^3 \bigg)   + 
\\
  +  \frac{\sigma}{4}\bigg(  \|s^{(i)} \|^4 +  4 \|s^{(i)} \|^2 s^Ts^{(i)}+ 2\|s^{(i)} \|^2 \|s \|^2 + 4  [s^T s^{(i)}]^2 + \|s\|^4 \bigg) + d_i\|s\|^4.
\end{eqnarray*}
\normalsize
Using $T[s^{(i)}+s]^3 - T[s]^3 =T[s^{(i)}]^3+ 3T[s]^2[s^{(i)}] + 3T[s^{(i)}]^2[s]   $ and $ \|s^{(i)}+s\|^4  =  \|s^{(i)}\|^4 + 4 \|s\|^2  s^Ts^{(i)} + 2 \|s\|^2 \|s^{(i)}\|^2 + 4[s^T s^{(i)}]^2 + 4 \|s^{(i)}\|^2  s^Ts^{(i)} +\|s\|^4$, we deduce that
\normalsize
\begin{eqnarray*}
M_c(s^{(i)}, s)  = f_0 + g^T[s^{(i)}+s]+\frac{1}{2}H [s^{(i)}+s]^2 + \frac{1}{6} \bigg( T[s^{(i)}+s]^3 -  T[s]^3  + \beta_i \|s\|^3 \bigg) +   \frac{ \sigma}{20} \big \|s^{(i)} + s \big\|^4   
\\
 + \frac{\sigma}{5} \bigg(\underbrace{ \|s^{(i)}\|^4 +  \underbrace{4 \|s^{(i)}\|^2  s^Ts^{(i)}  + 2 \|s\|^2 \|s^{(i)}\|^2}_{\ge 0 \text{ by \eqref{norm eq}}} + 4[s^T s^{(i)}]^2  + \|s \|^4- \|s\|^2  s^Ts^{(i)}  }_{I_4}\bigg)   + \underbrace{d_i  \|s\|^4}_{\ge 0}  . 
\end{eqnarray*}
Note that 
$$
I_4 \ge \|s^{(i)}\|^4  +4[s^T s^{(i)}]^2  + \|s \|^4 - \|s\|^2  s^Ts^{(i)}  \ge  \|s\|^4  - \|s\|^3 \|s^{(i)} \|   + \|s^{(i)} \|^4. 
$$
If $\|s\|\le\|s^{(i)} \| $, 
the sum of the last two terms is positive. If $\|s\|\ge\|s^{(i)} \| $ ,   the sum of the first two terms is positive.   In both cases, $I_4 \ge 0$, thus, we obtained the desired result. 
\end{proof}

Now, we are ready to prove that the iterates generated by Algorithm \ref{TCQR} are uniformly bounded.

\begin{theorem}
\label{thm bound on si}
Let $M_c$ be defined as in \eqref{t mc} with $d_i\ge 0$ and $|\beta_i| \le B$. 
Suppose that the CQR algorithmic framework (Algorithm \ref{TCQR}) or its variants are used with $s^{(0)} = 0$ and $s_c^{(i)} = \argmin_{s\in \R^n} M_c(s^{(i)}, s)$.   Then, for all  $i \ge 0$, 
\begin{eqnarray}
\big\| s^{(i)}+ s_c^{(i)}\big\| < r_c:= \max  \bigg\{ \bigg(\frac{60\|g\|}{\sigma}\bigg)^{1/3}, \bigg(\frac{30\|H\|}{\sigma}\bigg)^{1/2}  , \bigg(\frac{80B}{\sigma} + \frac{90\Lambda_0}{\sigma}\bigg) \bigg\}
\label{bound for si}
\end{eqnarray}
where $\|g\| = \sqrt{g^Tg}$, $\|H\| = \max_{i}|\lambda_{i}(H)|$ is the maximum absolute value of eigenvalues of $H$ and $\Lambda_0 := \|T\|=\max_{\|u_1\|=\|u_2\|=\|u_3\|=1} {T [u_1][u_2][u_3]}$. 
\end{theorem}

\begin{remark}
Note that $r_c$ is an iteration-independent bound and only depends on $g$, $H$, $T$, and $B$, which are fixed for $m_3$. The only requirement for Theorem \ref{thm bound on si} is that $|\beta_i| \le B$, where $B$ is a positive constant (i.e., the coefficients of the cubic-order term of $M_c$ are uniformly bounded), and $m_3(s^{(i)}) \le m_3(s^{(0)})$. The proof of Theorem \ref{thm bound on si} does not require the performance ratio test and is thus valid for both successful and unsuccessful iterations.
\end{remark}

\begin{proof}
    We prove the desired result by induction. For $i=0$, $s^{(0)} = 0$,  $M_c(s^{(0)}, 0) = f_0$, $s_c^{(0)} = \argmin_{s\in \R^n} M_c(s^{(0)}, s)$. 
Clearly, $0 \ge  M_c(s^{(0)}, s_c^{(0)}) - f_0$. Therefore, 
\begin{eqnarray*}
0 \ge g^T s_c^{(0)}+\frac{1}{2}H [s_c^{(0)}]^2 + \frac{\beta_0}{6}\|s_c^{(0)}\|^3   + \frac{\sigma}{4}\big \|s_c^{(0)}\big\|^4 \ge  -\|g\| \|s_c^{(0)}\| -\frac{\|H\|}{2}  \|s_c^{(0)}\|^2  - \frac{B}{6}\|s_c^{(0)}\|^3 + \frac{\sigma}{4}\big \|s_c^{(0)}\big\|^4
\end{eqnarray*}
where the second inequality uses the Cauchy-Schwarz inequality, norm properties, and $ \beta_0  > -B$. We further deduce that
\begin{eqnarray*}
0 \ge \bigg( -\|g\| +  \frac{\sigma}{12}\big \|s_c^{(0)}\big\|^3\bigg)\|s_c^{(0)}\| + \bigg(-\frac{\|H\|}{2}+  \frac{\sigma}{12}\big \|s_c^{(0)}\big\|^2\bigg)  \|s_c^{(0)}\|^2+ \bigg( - \frac{B}{6}+  \frac{\sigma}{12}\big \|s_c^{(0)}\big\|\bigg) \|s_c^{(0)}\|^3. 
\end{eqnarray*}
The inequality above cannot hold unless at least one of the terms on the right-hand side is negative, which is equivalent to 
$$
\big\|s^{(0)}+s^{(0)}_c\big\| \le \max \bigg\{ \bigg(\frac{12\|g\|}{\sigma}\bigg)^{1/3}, \bigg(\frac{6\|H\|}{\sigma}\bigg)^{1/2}  , \frac{2B}{\sigma}  \bigg\} < r_c. 
$$
Thus, \eqref{bound for si} is true for $i=0$. 

For the inductive hypothesis\footnote{This is the same as$\big\|s^{(i-1)}\| = \big\|s^{(i-1)} + s_c^{(i-1)}\big\| \le r_c$ if iteration is successful.}, assume \eqref{bound for si} is true at the $(i-1)$th iteration $\big\|s^{(i-1)} + s_c^{(i-1)}\big\| \le r_c$. 

At the $i$th iteration, either we are in the good case, $\big\|s^{(i)} + s_c^{(i)}\big\| \le \big\|s^{(i)}\big\|$ which leads to  $\big\|s^{(i)} + s_c^{(i)}\big\|  \le \big\|s^{(i)}\big\| \le r_c$ directly. Or, we are in the hard case, $\big\|s^{(i)} + s_c^{(i)}\big\| \ge \big\|s^{(i)}\big\|$. In this case, since $s_c^{(i)} = \argmin_{s \in \R^n} M_c(s^{(i)}, s)$, we have $0 \ge  M_c(s^{(i)},s_c^{(i)}) - M_c(s^{(i)},0) = M_c(s^{(i)},s_c^{(i)}) - f_i \ge \dotsc \ge M_c(s^{(i)},s_c^{(i)}) - f_0$. 
Denoting $\hat{s} := s^{(i)} + s_c^{(i)}$, using Lemma \ref{technical lemma}, we have
 \begin{eqnarray*}
  0 &\ge& M_c(s^{(i)}, s_c^{(i)}) -f_0 \ge  g^T \hat{s}+\frac{1}{2}H [\hat{s}]^2 + \frac{1}{6} T[\hat{s}]^3  -  \frac{1}{6}T[s_c^{(i)}]^3 + \frac{\beta_i}{6}\|s_c^{(i)}\|^3 + \frac{\sigma}{20}\big \|\hat{s} \big\| 
\\  &>&  -\|g\| \|\hat{s}\| -\frac{\|H\|}{2} \|\hat{s}\|^2 - \frac{ \Lambda_0}{6}\|\hat{s}\|^3  -  \frac{\Lambda_0} {6}\|s_c^{(i)}\|^3- \frac{B}{6}\|s_c^{(i)}\|^3  + \frac{\sigma}{20}\big \|\hat{s}\big\|^4. 
\end{eqnarray*}
The last inequality uses the Cauchy-Schwarz inequality, norm properties, and $\beta_i > -B.$
Since $\big\|\hat{s}\big\| \ge \big\|s^{(i)}\big\|, $
this gives 
$
     \|s_c^{(i)}\| =\big\|\hat{s} - s^{(i)}\big\|\le \big\|\hat{s}\big\| + \big\|s^{(i)}\big\| \le 2 \big\|\hat{s}\big\|.
$
Consequently, 
 \begin{eqnarray*}
  0 &>& -\|g\| \|\hat{s}\| -\frac{\|H\|}{2} \|\hat{s}\|^2 - \frac{3\Lambda_0}{2}\|\hat{s}\|^3   - \frac{4B}{3} \big \|\hat{s}\big\|^3 + \frac{\sigma}{20} \big \|\hat{s}\big\|^4
\\&=& \bigg( -\|g\| +  \frac{\sigma}{60}\big \|\hat{s}\big\|^3\bigg) \|\hat{s}\| + \bigg(-\frac{\|H\|}{2} +  \frac{\sigma}{60}\big \|\hat{s}\big\|^2\bigg)  \|\hat{s}\|^2  + \bigg( - \frac{4B}{3} - \frac{3\Lambda_0}{2}+  \frac{\sigma}{60}\big \|\hat{s}\big\|\bigg) \|\hat{s}\|^3. 
\end{eqnarray*}
The inequality above cannot hold unless at least one of the terms on the right-hand side is negative, which is equivalent to 
$
\big\|s^{(i)}+s_c^{(i)}\big\|= \|\hat{s}\big\| \le \max \bigg\{ \bigg(\frac{60\|g\|}{\sigma}\bigg)^{1/3}, \bigg(\frac{30\|H\|}{\sigma}\bigg)^{1/2}  , \bigg(\frac{80B}{\sigma} + \frac{90\Lambda_0}{\sigma}\bigg) \bigg\} = r_c. 
$
\end{proof}

\section{Proof of Lemma \ref{lemma Number of Iterations}}
\label{appendix: proof for lemma 3.6}
\begin{proof}
The update of regularization parameters in Algorithm \ref{algo: cqr variant 1}  gives that, for each $i$, 
$$
0\le d_{j+1} \le d_j, \quad j \in \mathcal{S}_i , \qquad \gamma \max\{1, d_j\} = d_{j+1}, \quad  j \in \mathcal{U}_i ,
$$
where $j \in [0:i]$. Using $d_0 = 0$, we inductively deduce that 
$
\gamma^{\mathcal|{U}_i |} = \max\{1, d_0\}  \gamma^{\mathcal|{U}_i |}\le d_i.
$ Therefore, using our assumption that $d_i \le d_{\max}$, we deduce that
$$
\mathcal|{U}_i | \log \gamma \le \log(d_{\max}). 
$$
The desired result follows from the inequality $i = \mathcal|{U}_i|+ \mathcal|{S}_i|$.
\end{proof}

\newpage
\section{Performance Profiles for Different $\beta$ Updates}
\label{appendix beta choice}
\begin{figure}[!h]
    \centering
        \caption{\small Performance profile of Algorithm \ref{PCQR} using the second choice of $\beta_i$. The same test set is used as in Figure \ref{fig beta choice 1}.}
    \label{fig beta choice 2}
    \includegraphics[width =15.5cm]{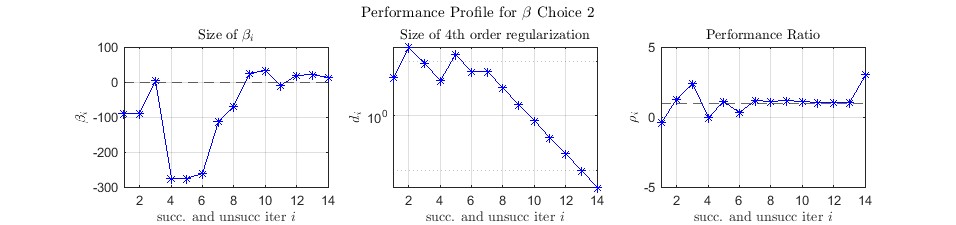}
\end{figure}

\begin{figure}[h]
    \centering
        \caption{\small Performance profile of Algorithm \ref{PCQR} using the third choice of $\beta_i$. The same test set is used as in Figure \ref{fig beta choice 1}.}
    \label{fig beta choice 3}
    \includegraphics[width =15.5cm]{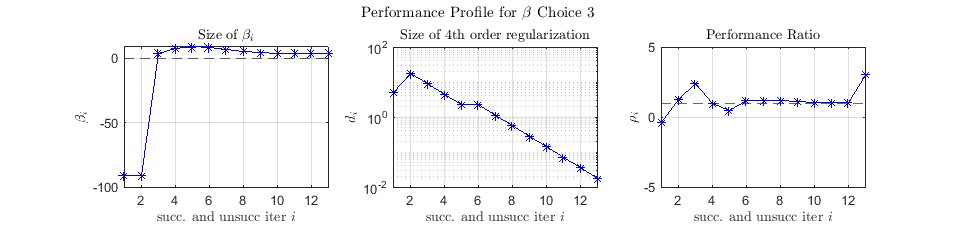}
\end{figure}

\section{Numerical Results for Various Test Sets}
\label{appendix: Numerical Results on Various Test Sets}

The results presented in each table are obtained over $10$ random problems randomized using MATLAB functions \texttt{rand()} or \texttt{randn()}. In all of the test cases, all algorithms converge to a point that satisfies $\|g_i\| < \epsilon$ with $\epsilon= 10^{-5}$.  The \textit{relative minimum values} of a method are obtained by dividing the minimum value achieved by this algorithm by the minimum of all the minimum values obtained by all algorithms. 

\begin{table}[!htbp]
\centering
\caption{\small \textbf{Convex $m_3$ and Convex $H$:} parameters in \eqref{numerical setu} are \texttt{a} $=80$, and $\sigma =80$.  In the \textbf{first part of the table}, $m_3$ is convex, $H = \texttt{symm(30*(randn(n)+neye(n)))}$ and \texttt{c} $=1$; In the \textbf{second part of the table}, $m_3$ is nonconvex but locally convex at $s=\boldsymbol{0}$, $H = \texttt{symm(30*(randn(n)+neye(n)))}$ and \texttt{c} $=80$. }  
    \label{table convex H}
{\small
\begin{tblr}{
  colspec = {c|rrr|rrr|rrr},
  cell{1}{2,9} = {c=3}{c}, 
}
\hline[2pt]
   & \SetCell[c=3]{}Relative Min. & &  &\SetCell[c=3]{}Iterations.\quad & & &  \SetCell[c=3]{}Evaluations.  \\
\hline[1pt]
Method  &   $n=5$ &  50 &  100 & 5 &  50 &  100  &   5 &  50 &  100 &\\
\hline[0.5pt]
    \textbf{\color{red}ARC}  &  1 & 1 & 1
    &  5.5&3&3
    &    5.5&3&3 &\\
    \textbf{Nesterov} &  1 & 1 & 1
    &  12.9 & 14 & 15 
    &  12.9 & 14 & 15 &\\
   \textbf{\color{blue}QQR-v1}  &  1 & 1 & 1
   &8.8 & 9 & 9
    &8.8 & 9 & 9 &\\
   \textbf{\color{blue}QQR-v2} &  1 & 1 & 1
    &  \bb{2.9} & \bb{2} & \bb{2} 
    &  \bb{2.9} & \bb{2} & \bb{2} &\\
   \textbf{\color{blue}CQR-choice 2}  &  1 & 1 & 1
    &  \bb{2.9} & \bb{2} & \bb{2} 
    &  \bb{2.9} & \bb{2} & \bb{2} &\\
   \textbf{\color{blue}CQR-choice 3}  &  1 & 1 & 1
    &  \bb{2.9} & \bb{2} & \bb{2} 
    &  \bb{2.9} & \bb{2} & \bb{2} &\\
\hline[0.5pt]
\hline[0.5pt]
    \textbf{\color{red}ARC}  &  1 & 1 & 1
    &  5.6 &4 &4
      &  5.6 &4 &4 &\\
   \textbf{\color{blue}QQR-v1}  &  1 & 1 & 1
   &8.8 & 9 & 9
    &8.8 & 9 & 9 \\
   \textbf{\color{blue}QQR-v2} &  1 & 1 & 1
    &  \bb{4} & \bb{3} & \bb{3} 
    &  \bb{4} & \bb{3} & \bb{3} & \\
   \textbf{\color{blue}CQR-choice 2}  &  1 & 1 & 1
    & \bb{4.9} & \bb{3} & \bb{3} 
    & \bb{4.7} & \bb{3} & \bb{3} & \\
   \textbf{\color{blue}CQR-choice 3}  &  1 & 1 & 1
    & \bb{4.7} & \bb{3} & \bb{3} 
    & \bb{4.7}  & \bb{3} & \bb{3} & \\
\hline[2pt]
\end{tblr}
}
\end{table}

\begin{table}[!htbp]
\centering
\caption{\small \textbf{Ill-Conditioned $H$, Diagonal $H$ or Singular $H$:} parameters in \eqref{numerical setu} are \texttt{a} $=80$, \texttt{c} $=80$ and $\sigma = 80$. $H$ is a diagonal ill-conditioned matrix with diagonal entries uniformly distributed in $[-10^{-10}, 10^{10}]$.   }
\label{table ill con H}
{\small
\begin{tblr}{
  colspec = {c|rrr|rrr|rrr},
  cell{1}{2,9} = {c=3}{c}, 
}
\hline[2pt]
   & \SetCell[c=3]{}Relative Min. & &  &\SetCell[c=3]{}Iterations.\quad & & &  \SetCell[c=3]{}Evaluations.  \\
\hline[1pt]
Method  &   $n=5$ &  50 &  100 & 5 &  50 &  100  &   5 &  50 &  100 & \\
\hline[0.5pt]
    \textbf{\color{red}ARC}  &  1 & 1 & 0.99033 
    &  9.6 & 8.6 & 9.3
    &  6.5 & 5.7&6& \\
   \textbf{\color{blue}QQR-v1}   &  1 & 1 & 0.99033  
   &9 & 9.7& 10.2
   &8.8 & 8.6 & 9.2 \\
   \textbf{\color{blue}QQR-v2}   &  1 & 1 & 0.99033  
   &\bb{4.1} & \bb{5} & \bb{4.7}
   &\bb{4.1} & \bb{4.9} & \bb{4.5} \\
   \textbf{\color{blue}CQR-choice 2}   &  1 & 1 & 1
   &\bb{4.9} & \bb{7.5} & \bb{6.5}
   &\bb{4.4} & \bb{5.3} & \bb{4.7} \\
   \textbf{\color{blue}CQR-choice 3}   &  1 & 1 & 0.99033 
   &\bb{4.7} & \bb{5.6} & \bb{5.7}
   &\bb{4.1} & \bb{4.7} & \bb{4.5} \\
\hline[2pt]
\end{tblr}
}
\end{table}

\begin{table}[!htbp]
\centering
    \caption{\small \textbf{Changing $\sigma$:}  From top to bottom, parameters in \eqref{numerical setu} are \texttt{a} $=80$, \texttt{b} $=80$,  \texttt{c} $=80$ and $\sigma =5, 300$, respectively.  The experiments show the algorithm's ability to handle a range of $\sigma$ values. Generally, smaller values of $\sigma$ for $m_3(s)$ require more function/derivative evaluations for convergence across all algorithms. }
     \label{change sigma}
{\small
\begin{tblr}{
  colspec = {c|rrr|rrr|rrr},
  cell{1}{2,9} = {c=3}{c}, 
}
\hline[2pt]
   & \SetCell[c=3]{}Relative Min. & &  &\SetCell[c=3]{}Iterations.\quad & & &  \SetCell[c=3]{}Evaluations.  \\
\hline[1pt]
Method  &   $n=5$ &  50 &  100 & 5 &  50 &  100  &   5 &  50 &  100 &\\
\hline[0.5pt]
    \textbf{\color{red}ARC}  &  1 & 1 & 0.98313
    &  11.8& 21.7 &26.3
    &  8.1 &15.1 &20.5 &\\
   \textbf{\color{blue}QQR-v1}  & 1 & 0.9958 & 1
   &14.5& 19.6 & 22
   &12.5 & 17.3 &19.5 &\\
   \textbf{\color{blue}QQR-v2} & 1 &  0.9958 &   0.97191
    &  \bb{8.4} & \bb{14} & \bb{16.5} 
    &  \bb{7.6} & \bb{12.1} & \bb{14.4} &\\
   \textbf{\color{blue}CQR-Choice 2}  & 1 & 1 & 0.98498 
&\bb{8.2}& \bb{15.8} & \bb{17.5}
&\bb{7} & \bb{12} &\bb{15.5} &\\
\textbf{\color{blue}CQR-Choice 3} & 1 & 1  & 0.97191  
&  \bb{8.2} & \bb{13.1} &\bb{18.8}
&  \bb{7.2} & \bb{11.8} & \bb{14.8} &\\
\hline[0.5pt]
\hline[0.5pt]
    \textbf{\color{red}ARC}  &  1 & 1 &  0.99642 
    & 12.5 & 18.4 &19.2
    &6.5 &8.5 & 12.2 &\\
   \textbf{\color{blue}QQR-v1}   &  1 & 1 & 0.99825
   &9.9& 12.6& 16.5 &
   8.5& 11.8& 15.3&\\
   \textbf{\color{blue}QQR-v2} &  1 & 1 & 1
    &\bb{4.9} & \bb{7.6}& \bb{10}
    &  \bb{4.9} & \bb{7.4} & \bb{9.7} &\\
       \textbf{\color{blue}CQR-Choice 2}  & 1 & 1 &  0.99642
   &\bb{4.9} & \bb{10.7}& \bb{11.2}
    &  \bb{4.9} & \bb{7.1} & \bb{9.6} &\\
\textbf{\color{blue}CQR-Choice 3} & 1 & 1  & 0.99642
   &\bb{4.8} & \bb{7}& \bb{16}
    &  \bb{4.8} & \bb{7} & \bb{9.8} &\\
\hline[2pt]
\end{tblr}
}

\end{table}

\begin{table}[!htbp]
\centering
    \caption{\small \textbf{Large Tensor Term:}  parameters in \eqref{numerical setu} are \texttt{a} $=80$, \texttt{b} $=80$, $\sigma =5$, and \texttt{c} $=300$ respectively. Minimizing $m_3(s)$ becomes progressively harder with a larger tensor term. }
        \label{table small tensor}
{\small
\begin{tblr}{
  colspec = {c|rrr|rrr|rrr},
  cell{1}{2,9} = {c=3}{c}, 
}
\hline[2pt]
   & \SetCell[c=3]{}Relative Min. & &  &\SetCell[c=3]{}Iterations.\quad & & &  \SetCell[c=3]{}Evaluations.  \\
\hline[1pt]
Method  &   $n=5$ &  50 &  100 & 5 &  50 &  100  &   5 &  50 &  100 & \\
\hline[0.5pt]
    \textbf{\color{red}ARC}  &  1 & 1& 0.9845
    &  16.7 & 24.3 & 35.1 
    &  9.1 & 15.8& 27.1 & \\
   \textbf{\color{blue}QQR-v1}   &   0.99204 &  0.94926 & 1
   &18.3& 22.1 & 25.6
   &15.2 & 18.6& 22.4 \\
   \textbf{\color{blue}QQR-v2}   &  1 & 1& 0.96128 
    &  \bb{10.6}& \bb{18.2}& \bb{22.1} 
    &  \bb{8.4} & \bb{12.6} & \bb{18} & \\
   \textbf{\color{blue}CQR-Choice 2}  &  1 & 1& 0.98301 
    &  \bb{10.8} & \bb{15.4} & \bb{22} 
    &  \bb{8.3} & \bb{13.1} & \bb{19.1}& \\
\textbf{\color{blue}CQR-Choice 3} &  1 & 1& 0.98301 
    &  \bb{11.2}& \bb{20.9} & \bb{19.6} 
    &  \bb{9} & \bb{13.2} & \bb{17.3} & \\
\hline[2pt]
\end{tblr}
}

\end{table}
 
\begin{table}[!htbp]
\centering
\caption{\small \textbf{Directional and Diagonal Tensors:} In the \textbf{first table}, $g$ and $T$ are negative and directional and $H$ is near zero where \texttt{g = -80*rand(n,1)} and \texttt{T = -80*symm(rand(n,n,n))}. The parameters in \eqref{numerical setu} are \texttt{b} $=0.1$ and $\sigma =80$. In the \textbf{second table}, we consider positive and directional tensors, where parameters in \eqref{numerical setu} are \texttt{a} $=80$, \texttt{b} $=80$, $\sigma =80$ and $T$ is a diagonal tensor with entries that are uniformly distributed in the range of $[0,40]$. }
     \label{table large tensor}
{\small
\begin{tblr}{
  colspec = {c|rrr|rrr|rrr},
  cell{1}{2,9} = {c=3}{c}, 
}
\hline[2pt]
   & \SetCell[c=3]{}Relative Min. & &  &\SetCell[c=3]{}Iterations.& & &  \SetCell[c=3]{}Evaluations.  \\
\hline[1pt]
Method  &   $n=5$ &  50 &  100 & 5 &  50 &  100  &   5 &  50 &  100 & \\
\hline[0.5pt]
    \textbf{\color{red}ARC}  &  1 & 1 & 1
    &  13.2 & 17.4 & 15.4
    &  7.1 & 8.9 & 7.7& \\
   \textbf{\color{blue}QQR-v1} &  1 & 1 & 1
   &11.7&19.3 & 21.1
   &10.5 & 15.3 & 16.6 \\
   \textbf{\color{blue}QQR-v2} &  1 & 1 & 1
    &  \bb{6.2} & \bb{8.7}  & \bb{12.1} 
    &  \bb{6.1} & \bb{7.3} & \bb{8.2} & \\
       \textbf{\color{blue}CQR-Choice 2}  &  1 & 1 & 1
    &  \bb{10.4} & \bb{8.7}  & \bb{8} 
   &  \bb{5.7} & \bb{7.7} & \bb{8} & \\
\textbf{\color{blue}CQR-Choice 3} &  1 & 1 & 1
    &  \bb{6.7} & \bb{8.4}  & \bb{11.2} 
   &  \bb{6.1} & \bb{7.4} & \bb{8} & \\
\hline[0.5pt]
\hline[0.5pt]
    \textbf{\color{red}ARC}  &  1 & 1 & 1 
    &  11.1 &12.5 &16.4
    &  6.2 & 6.5 & 8 & \\
   \textbf{\color{blue}QQR-v1}   &  1 & 1 & 1 
   &10.1  & 10.5 & 10.4
   &8.7 & 8.9 & 8.4 \\
   \textbf{\color{blue}QQR-v2} &  1 & 1 & 1 
    &\bb{4.1} & \bb{4.4}& \bb{4.1}
    &  \bb{4.1} & \bb{4.4} & \bb{4.1}& \\
       \textbf{\color{blue}CQR-Choice 2}  & 1 & 1 & 1 
    &\bb{8.3} & \bb{5.4}& \bb{5.2}
    &  \bb{4.9} & \bb{5.4} & \bb{5.2}& \\
\textbf{\color{blue}CQR-Choice 3} &  1 & 1 & 1 
    &\bb{4.9} & \bb{5.4}& \bb{5.2}
    &  \bb{4.9} & \bb{5.4} & \bb{5.2}& \\
\hline[2pt]
\end{tblr}
}
\end{table}

\clearpage
\newpage
\section{More Examples of Badly Scaled Subproblems}
\label{appendix: Special Examples}

\begin{figure}[!ht]
    \centering
\caption{\small \textbf{Examples with varying magnitudes of small entries:} In these tests, the setup is the same as in Figure \ref{bad scale 0.1} with $n=99$. The left two plots depict scenarios where the other entries follow a standard normal distribution $\mathcal{N}(0, 1)$. In the right two plots, all other entries, except the first one, are set to zero.}
    \label{bad scale 1}
        \includegraphics[width =7.5cm]{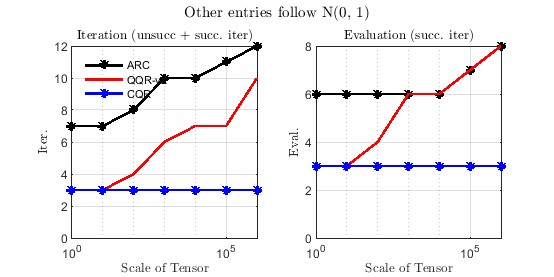}
    \includegraphics[width =7.5cm]{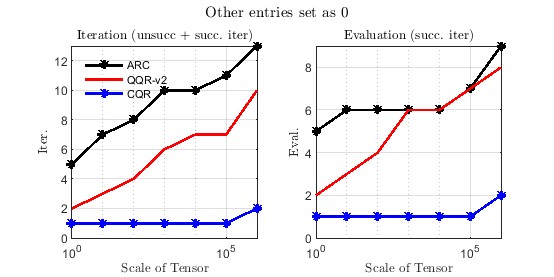}
         \label{fig special example 0}
\end{figure}

\begin{figure}[!ht]
    \centering
\caption{\small \textbf{Examples with different $\beta$ Choices:} These scenarios use the same setup as in Figure \ref{bad scale 0.1}, but $\beta_i$ is chosen according to choice 3. Similar performance advantages are observed.}
    \label{choice 3}
    \includegraphics[width =7.5cm]{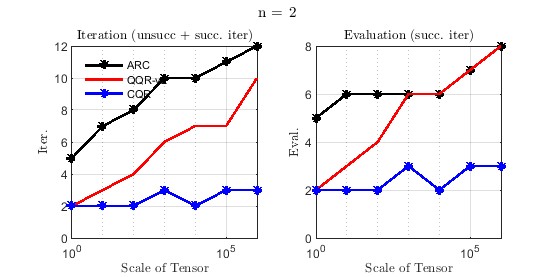}
    \includegraphics[width =7.5cm]{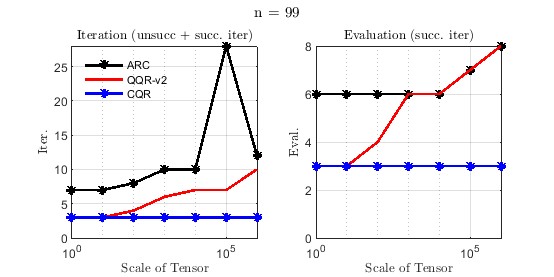}
\end{figure}

\begin{figure}[!ht]
    \centering
\caption{\small \textbf{Examples with multiple large entries:} These tests maintain the same setup as in Figure \ref{bad scale 0.1}, with $n=99$. In the left two plots, apart from the first entry diagonal of $T$, the size of certain off-diagonal entries, such as $T[2,3,4]$ and its permutations, also varies from $-0.1$ to $-10^5$. In the right two plots, aside from the first diagonal entry of $T$, the size of two diagonal entries of $T$ is also varied from $-0.1$ to $-10^5$.}
        \includegraphics[width =7.5cm]{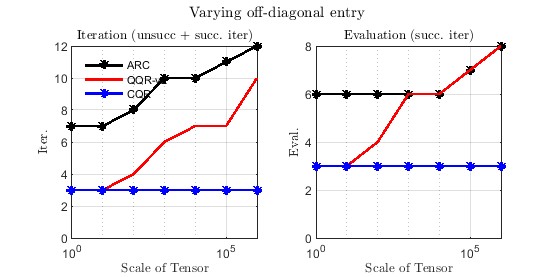}
    \includegraphics[width =7.5cm]{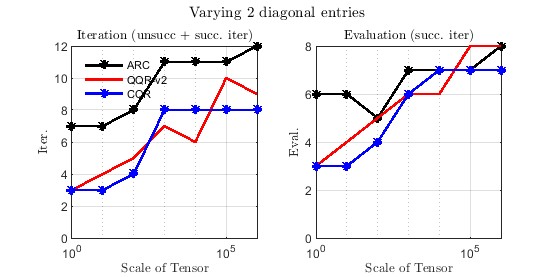}
\label{fig special example diag}
\end{figure}


\scriptsize{
\bibliographystyle{plain}
\bibliography{sample.bib}
}
\end{document}